\documentclass[
final
]{dmtcs-episciences}

\usepackage[utf8]{inputenc}
\usepackage{subfigure}

\usepackage[round]{natbib}
\usepackage{amsfonts}
\usepackage{amsmath}
\usepackage{mathrsfs}
\usepackage{mathrsfs,amscd,amssymb,color,amsthm,amsmath,bm,graphicx,psfrag,subfigure,url}
\newtheorem{thm}{Theorem}[section]
\newtheorem{lem}[thm]{Lemma}
\newtheorem{cor}[thm]{Corollary}

\newenvironment{wst}
{\setlength{\leftmargini}{1.5\parindent}
 \begin{itemize}
 \setlength{\itemsep}{-1.1mm}}
{\end{itemize}}
\author{Huihui Zhang\affiliationmark{1,3}\thanks{She is fully supported by the excellent doctoral  dissertation  cultivation grant from Central China Normal University (Grant No. 2016YBZZ084).}
  \and Jing Chen\affiliationmark{2}\thanks{She is fully supported by the National Natural Science Foundation of China (Grant No.11501188).}
  \and Shuchao Li\affiliationmark{3}\thanks{Corresponding author and I am fully supported by the National Natural Science Foundation of China (Grant No. 11671164).}}
\title[The (revised) Szeged index and Wiener index]{On the quotients between the (revised) Szeged index and Wiener index of
graphs}
\affiliation{
  Department of Mathematics, Luoyang Normal Univeristy, Luoyang 471002, P.R. China\\
  Department of Mathematics, Hunan First Normal University, Changsha, Hunan 410205, P.R. China\\
  Faculty of Mathematics and Statistics,  Central China Normal
University, Wuhan 430079, P.R. China}
\keywords{Szeged index; Revised Szeged index; Wiener index}
\received{2016-6-26}
\revised{2017-3-6}
\accepted{2017-4-14}
\begin{document}
\publicationdetails{19}{2017}{1}{12}{1514}
\maketitle

\begin{abstract}
  Let $Sz(G),Sz^*(G)$ and $W(G)$ be the Szeged index, revised Szeged index and Wiener index of a
graph $G.$ In this paper, the graphs with the fourth, fifth, sixth and seventh largest Wiener indices among all unicyclic graphs of order $n\geqslant 10$ are characterized; and the graphs with the first, second, third, and fourth largest Wiener indices among all bicyclic graphs are identified. Based on these results, further relation on the quotients between the (revised) Szeged index and the Wiener index are studied. Sharp lower bound on $Sz(G)/W(G)$ is determined for all connected graphs each of which contains at least one non-complete block. In addition the connected graph with the second smallest value on $Sz^*(G)/W(G)$ is identified for $G$ containing at least one cycle.
\end{abstract}

\section{Introduction}
\label{sec:in}
We consider that all graphs in this paper are finite, undirected and simple. We follow the notations and terminologies in \cite{3} except otherwise stated. Let $G=(V_G, E_G)$ be a connected graph with vertex set $V_G$ and edge set $E_G.$ A connected graph is \textit{cyclic} if it contains at least one cycle. In particular, a connected graph $G$ is \textit{unicyclic} (resp. \textit{bicyclic}) if $|E_G|=|V_G|$ (resp. $|E_G|=|V_G|+1).$ For convenience, let $|G|:=|V_G|.$

In the subsequent sections, we use $G-v$, or $G-uv$ to denote the graph obtained from $G$ by deleting vertex $v \in V_G$, or edge
$uv \in E_G$, respectively (it is naturally extended if at least two vertices or edges are deleted). Let
$G+uv$ be the graph obtained from $G$ by adding an edge $uv \not\in E_G$. For a subset $S$ of $V_G,$ let $G[S]$ be the subgraph induced by $S$. For $v\in V_G$, we denote by $N_G(v)$ (or $N(v)$ for short) the set of all neighbors of $v$ in $G$ and let $d_G(v)=|N_G(v)|$ be the degree of $v$ in $G.$ Call $u$ a \textit{pendant vertex} or \textit{leaf} in $G$, if $d_G(u)=1.$ We denote by $P_n,C_n, S_n$ and $K_n$ the path, cycle, star and complete graph of order $n$, respectively. We call $L_{n,r}$ a \textit{lollipop} if it is obtained by identifying some vertex of $C_r$ with an end-vertex of $P_{n-r+1}.$

For $u, v\in V_G$, the \textit{distance} between $u$ and $v$ in $G$, denoted by $d_G(u,v),$ is the length of a shortest path connecting $u$ and $v$. The \textit{diameter}, $d(G),$ of $G$ is equal to $\max_{u,v\in V_G}d_G(u, v).$ For all $v\in V_G$, let $\eta_G^i(v)=|\{u\in V_G: d(u,v)=i\}|.$ The symbol $\cong$ denotes that two graphs in question are isomorphic.

The \textit{Wiener index} of $G$ is defined as the sum of all distances between pairs of unordered vertices in $G$, i.e.,
\begin{equation}
 W(G)=\sum_{\{u,v\}\subseteq V_G}d_G(u, v)=\frac{1}{2}\sum_{u\in V_G}D_G(u), \label{eq:1.1}
\end{equation}
where $D_G(u)=\sum_{x\in V_G}d_G(x, u).$ This distance-based graph invariant was in chemistry introduced back in \cite{25} and in mathematics about 30 years later;~\cite{008}. Nowadays, the Wiener index is a extensively studied graph invariant; see the reviews~\cite{D-E-G, D-G-K-Z}. A collection of recent papers dedicated to the investigations of the Wiener index; see~\cite{KM-S,LS-S,LR-D}. 
The problem of finding an upper bound on the Wiener index of a graph is quite challenging; see \cite{M-V}. Only a few papers considered the upper bounds on the Wiener index of graphs; see \cite{D,Z,T-D,D-Z}.

Given an edge $e=uv$ in $G$, define three sets with respect to $e$ as follows:
$$
\begin{array}{cc}
  N_u(e)=\{w\in V_G:d_G(u, w)< d_G(v, w)\}, & N_v(e)=\{w\in V_G:d_G(v, w)< d_G(u, w)\}, \\[5pt]
  N_0(e)=\{w\in V_G:d_G(u, w)= d_G(v, w)\}. &
\end{array}
$$
Clearly, $V_G=N_u(e)\cup N_v(e)\cup N_0(e)$. For convenience, let $n_u(e)=|N_u(e)|,\, n_v(e)=|N_v(e)|$ and $n_0(e)=|N_0(e)|$. It is easy to see $n_u(e)+n_v(e)+n_0(e)=|V_G|$. If $G$ is bipartite, then $N_0(e)=\emptyset$ holds for all $e\in E_G$. Consequently, for any bipartite graph $G$ with $e\in E_G$, $n_u(e)+n_v(e)=|V_G|$.

From \cite{D-E-G,12,25}, we know that, for a tree $T$, its Wiener index can be defined alternatively as
\begin{equation}\label{eq:1.2}
  W(T)=\sum_{e=uv\in E_T}n_u(e)n_v(e).
\end{equation}
Motivated by (\ref{eq:1.2}), \cite{10} introduced the \textit{Szeged index} of graph $G$, which is defined by
\[\label{eq:1.3}
 Sz(G)=\sum_{e=uv\in E_G}n_u(e)n_v(e).
\]

\cite{23} observed that the vertices at equal distances from the end-vertices of an edge do not contribute to the Szeged index, and so he proposed the \textit{revised Szeged index} $Sz^*(G)$ of a  graph $G$ as follows:
$$
  Sz^*(G)=\sum_{e=uv\in E_G}\left(n_u(e)+\frac{n_0(e)}{2}\right)\left(n_v(e)+\frac{n_0(e)}{2}\right).
$$
For more recent results on (revised) Szeged index, one may be referred to these in \cite{2,15,28,18,21,24,26}.

By (\ref{eq:1.2}) and (\ref{eq:1.3}), we know that $Sz(T)= W(T)$ holds for any tree $T$. Then, many researchers focused on the difference between the Szeged index and the Wiener index on general graphs. Given a graph $G$, the difference $Sz(G)-W(G)\geqslant 0$ holds, which was conjectured in \cite{10} and proved in \cite{klavzar-1996}. Moreover, $Sz(G)=W(G)$ holds if and only if every block of $G$ is a complete graph, which was obtained by \cite{7}, and see \cite{16} for another proof. \cite{19} studied the structure of graphs $G$ with $Sz(G)-W(G)=k$, here $k$ is a positive integer. In particular, \cite{20} identified the graphs for which the difference is 4 and 5. The difference between $Sz(G)$ and $W(G)$ in networks was investigated in~\cite{17}. \cite{21} showed that, if $G$ is connected, then $Sz^*(G)\geqslant Sz(G)$ with equality if and only if $G$ is bipartite. Some further results on the difference  between the Wiener index and the (revised) Szeged index were established in~\cite{Z-L-Z}.

The computer program AutoGraphiX was used to study the relationship involving graph invariants; see \cite{1,4,9} for more detailed information. \cite{Han} used the computer program AutoGraphiX to generate eight conjectures on the difference (resp. quotient) between the (revised) Szeged index and Wiener index. \cite{5,C-L-L} confirmed three conjectures on the difference between the (revised) Szeged index and Wiener index, which can be summarized as following three theorems.
\begin{thm}[\cite{5,C-L-L}]\label{thm1.1}
Let $G$ be a connected bipartite graph with $n\geqslant 4$ vertices and $|E_G|\geqslant n$ edges. Then
$Sz(G)-W(G)\geqslant 4n-8.$
The equality holds if and only if $G$ is composed of a cycle $C_4$ on $4$ vertices and a tree $T$ on $n-3$ vertices sharing a single vertex.
\end{thm}
\begin{thm}[\cite{C-L-L}]\label{thm1.2}
Let $G$ be a connected graph with $n\geqslant 5$ vertices with an odd cycle and girth $g\geqslant 5$. Then
$Sz(G)-W(G)\geqslant 2n-5.$
The equality holds if and only if $G$ is composed of a cycle $C_5$ on 5 vertices, and one tree rooted at a vertex of the $C_5$ or two trees, respectively, rooted at two adjacent vertices of the $C_5.$
\end{thm}
\begin{thm}[\cite{C-L-L}]\label{thm1.3}
Let $G$ be a connected graph with $n \geqslant 4$ vertices and $|E_G|\geqslant n$ edges and with an odd cycle.
Then
$$
  Sz^*(G)-W(G)\geqslant \frac{n^2+4n-6}{4}.
$$
The equality holds if and only if $G$ is composed of a cycle $C_3$ on $3$ vertices and a tree $T$
on $n-2$ vertices sharing a single vertex.
\end{thm}
Recently, \cite{L-Z} confirmed three additional above conjectures, which are described as the following three theorems.
\begin{thm}[\cite{L-Z}]\label{lem:2.4}
Let $G$ be a cyclic graph of order $n\geqslant 4.$ 
\begin{wst}
\item[{\rm (i)}] If $G$ is a bipartite graph, then
$$
\frac{Sz^*(G)}{W(G)}\geqslant 1+\frac{24(n-2)}{n^3-13n+36}
$$
with equality if and only if $G$ is the lollipop $L_{n,4}.$
\item[{\rm (ii)}] If $G$ is a non-bipartite graph, then
$$
\frac{Sz^*(G)}{W(G)}\geqslant 1+\frac{3(n^2+4n-6)}{2(n^3-7n+12)}
$$
with equality if and only if $G$ is the lollipop $L_{n,3}.$
\end{wst}
\end{thm}
\begin{thm}[\cite{L-Z}]
Let $G$ be a unicyclic graph on $n\geqslant 4$ vertices. Then
$$\frac{Sz(G)}{W(G)}\leqslant\left\{
  \begin{array}{ll}
    2-\frac{8}{n^2+7}, & \hbox{if $n$ is odd,} \\
    2, & \hbox{if $n$ is even}
  \end{array}
\right.$$
with equality if and only if $G$ is the lollipop $L_{n,n-1}$ if $n$ is odd and the cycle $C_n$ if $n$ is even.
\end{thm}
\begin{thm}[\cite{L-Z}]
Let $G$ be a unicyclic graph on $n\geqslant 4$ vertices. Then
$$\frac{Sz^*(G)}{W(G)}\leqslant\left\{
  \begin{array}{ll}
    2+\frac{2}{n^2-1}, & \hbox{if $n$ is odd,} \\
    2, & \hbox{if $n$ is even}
  \end{array}
\right.$$
with equality if and only if $G$ is the cycle $C_n.$
\end{thm}

\cite{L-Z} determined sharp lower bounds on $Sz(G)/W(G)$ for cyclic graph with girth at least $4$.
\begin{thm}[\cite{L-Z}]\label{lem:2.3}
Let $G$ be a cyclic graph of order $n\geqslant 5$ with girth at least $4.$
\begin{wst}
\item[{\rm (i)}] If $G$ is a bipartite graph, then
$$
\frac{Sz(G)}{W(G)}\geqslant 1+\frac{24(n-2)}{n^3-13n+36}
$$
with equality if and only if $G$ is the lollipop $L_{n,4}.$
\item[{\rm (ii)}] If $G$ is a non-bipartite graph, then
$$
\frac{Sz(G)}{W(G)}\geqslant 1+\frac{6(2n-5)}{n^3-25n+90}
$$
with equality if and only if $G$ is the lollipop $L_{n,5}.$
\end{wst}
\end{thm}

\cite{8} showed that $Sz(G)/W(G)\geqslant 1$ with equality if and only if every block of $G$ is  complete, whereas the case that $G$ contains at least one block being non-complete is open. Based on Theorem~1.4, we may induce that $L_{n,4}$ is the unique graph with the smallest value on $Sz^*(G)/W(G)$ among $n$-vertex cyclic graphs. How can we determine the graph with second smallest value on $Sz^*(G)/W(G)$ among $n$-vertex cyclic graphs? Theorem~1.7 only considers the case for the graph with girth at least 4, whereas the case for cyclic graph of girth~3 is still open. In order to give solutions for the above open problems, it is natural and interesting for us to study some further relation on the quotients between the (revised) Szeged index and Wiener index of connected graphs.

Our paper is organized as follows: In Section 2, we give necessary definitions and state the main results of the paper. The first result determines the graphs with the fourth, fifth, sixth and seventh largest values on $W(G)$ among all unicyclic graphs of order $n\geqslant 10;$ while the second result characterizes the graphs with the first four greatest values on $W(G)$ in the class of all bicyclic graphs. These two results will be proved in Sections 4 and 5, respectively.
The third result obtains the sharp lower bound on $Sz(G)/W(G)$ for all connected graphs each of which contains at least one non-complete block. It, respectively, extends the results obtained by \cite{8} and \cite{L-Z}. The last result determines the cyclic graph $G$ with the second smallest value on $Sz^*(G)/W(G)$, which extends the result obtained by \cite{L-Z}. These two theorems are then proved in Section 6. In Section 3, we give some preliminary results which are used to prove our main results.
\section{Main results}
Consider a cycle $C_r$ whose vertices are labeled consecutively by $v_1,v_2,\ldots,v_r.$ Then let $C_r(T_1,T_2,\ldots,T_r)$ be an $n$-vertex graph obtained from $C_r$ and rooted trees $T_i$'s by identifying the root, say $r_i$, of $T_i$ with $v_i$ on $C_r, i=1,2,\ldots, r.$ Assume that $|V_{T_i}\setminus \{r_i\}|=n_i.$ Clearly, for $i=1,2,\ldots, r,\, n_i\geqslant 0.$ Thus, $|C_r(T_1,T_2,\ldots,T_r)|=\sum_{i=1}^r n_i+r.$ In particular, if every rooted tree is a path whose root is just one of its pendant vertices, then we denote $C_r(T_1,T_2,\ldots,T_r)$ by $C_r(P_{n_1+1},P_{n_2+1},\ldots,P_{n_r+1}).$
\begin{figure}[h!]
\begin{center}
\psfrag{a}{$H_n^0$}\psfrag{b}{$H_n^1$}\psfrag{C}{$H_{n}^3$}\psfrag{d}{$H_n^2$}
\includegraphics[width=120mm]{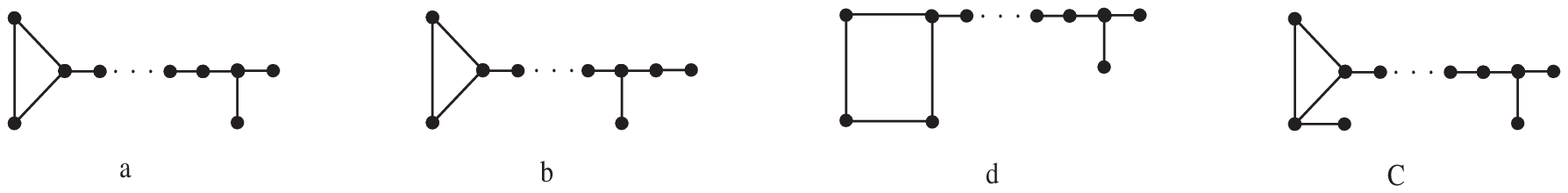}
  \caption{Graphs $H_n^0, H_n^1, H_n^2$ and $H_{n}^3$.}\label{fig2}
\end{center}
\end{figure}
For simplicity, let $H_n^0, H_n^1, H_n^2$ and $H_n^3$ be four unicyclic graphs of order $n$ as depicted in Fig. \ref{fig2}. Then, let $\mathscr{A}_n=\{L_{n,3}, L_{n,4}, H_n^0,H_n^1,H_n^2, H_n^3, C_3(P_{n-3},P_2,P_1),C_3(P_{n-4},P_3,P_1), C_4(P_{n-4},P_1,P_2,P_1)\}.$

\cite{T-D} showed that, among the unicyclic graphs of order $n\geqslant 6$, the graph $L_{n,3}$ is the graph with the largest Wiener index, both $C_3(P_{n-3},P_2,P_1)$ and $L_{n,4}$ are the graphs with the second largest Wiener index and $H_n^0$ is the graph with the third largest Wiener index.
Our first main result in this paper characterizes the unicyclic graphs of order $n\geqslant 10$ with the fourth, fifth, sixth and seventh largest Wiener indices.
\begin{thm}\label{thm5.1}
Let $G$ be a unicyclic graph of order $n\geqslant 10$ and $G$ is not in $\mathscr{A}_{n}.$
\begin{wst}
\item[{\rm (i)}]If $n=10$ and $G\ncong C_3(P_5,P_4, P_1),$ then
$W(G)<W(H_{10}^3)=W(H_{10}^2))<W(H_{10}^1)=W(C_3(P_5,\linebreak P_4,P_1))<W(C_4(P_6,P_1,P_2,P_1))<W(C_3(P_6, P_3,P_1)).$
\item[{\rm (ii)}]If $n=11$ and $G\ncong C_3(P_6,P_4, P_1),$ then $W(G)<W(H_{11}^3)=W(H_{11}^2)=W(C_3(P_6,P_4, P_1))<W(H_{11}^1) <W(C_4(P_7,P_1, P_2,P_1)) <W(C_3(P_7, P_3,P_1)).$
\item[{\rm (iii)}]If $n\geqslant  12,$ then $W(G)<W(H_n^3)=W(H_n^2)<W(H_n^1)<W(C_4(P_{n-4},P_1,P_2,P_1))<\linebreak W(C_3(P_{n-4},P_3,P_1)).$
\end{wst}
\end{thm}

Let $K_4^-$ be the graph obtained from $K_4$ by deleting one of its edges. Then, let $B_n^{(1)}$
be the $n$-vertex graph obtained by identifying an end-vertex of $P_{n-3}$ with a 3-degree vertex in $K_4^-$. For $0 \leqslant s \leqslant \lfloor\frac{n-4}{2}\rfloor,$ let $B(n,s)$ be the $n$-vertex graph obtained by attaching paths $P_{n-s-4}$ and $P_s$, respectively, to two 2-degree vertices of $K_4^-.$ Given two cycles $C_p, C_q,$ let $B_n^{p,q}$ be the $n$-vertex graph obtained by joining $C_p$ and $C_q$ with a path of length $n-p-q+1$.
Our next main result characterizes all the $n$-vertex bicyclic graphs with the first four largest Wiener indices.
\begin{thm}\label{thm5.2}
Let $G$ be a bicyclic graph of order $n\geqslant 6$ and $G$ is not in $\{B(n,0),B(n,1),B_n^{(1)},B_n^{3,3}\}.$
\begin{wst}
\item[{\rm (i)}]If $n=6,$ then $W(G)<W(B_6^{(1)})<W(B(6,1))=W(B_8^{3,3})<W(B(6,0)).$
\item[{\rm (ii)}]If $n=8$ and $G\ncong B(8,2),$ then $W(G)<W(B(8,2))=W(B_8^{(1)})<W(B(8,1))<W(B_8^{3,3})\linebreak <W(B(8,0)).$
\item[{\rm (iii)}]If $n\geqslant7$ and $n\not=8,$ then $W(G)<W(B_n^{(1)})<W(B(n,1))<W(B_n^{3,3})<W(B(n,0)).$
\end{wst}
\end{thm}

Based on Theorem 2.2, our next main result determines the sharp lower bound on $Sz(G)/W(G)$ for all connected graphs each of which contains at least one non-complete block.
\begin{thm}\label{thm4.1}
Let $G$ be an $n$-vertex connected graph containing at least one non-complete block, $n\geqslant 5$. Then
$$
\frac{Sz(G)}{W(G)}\geqslant 1+\frac{12}{n^3-19n+54}
$$
with equality if and only if $G\cong B_n^{(1)}.$
\end{thm}

Based on Theorem 2.1, our last main result characterizes the connected cyclic graph $G$ with the second smallest value on $Sz^*(G)/W(G)$.
\begin{thm}\label{thm4.2}
Let $G\ (\ncong L_{n,4})$ be a cyclic graph on $n\geqslant 10$ vertices.
Then
$$
\frac{Sz^*(G)}{W(G)}\geqslant 1+\frac{24(n-2)}{n^3-19n+54}
$$
with equality if and only if $G\cong H_n^2,$ where $H_n^2$ is depicted in Fig. \ref{fig2}.
\end{thm}

\section{Preliminaries}
In this section, we give some preliminary results and definitions which are used to prove our main results in the subsequent sections.

\begin{lem}[\cite{X-D}]\label{lem:2.6}
Let $G$ be an $n$-vertex bicyclic graph of diameter $n-2.$ Then
$G\in \{B(n,s): 0 \leqslant s \leqslant \lfloor\frac{n-4}{2}\rfloor\}.$
\end{lem}
\begin{figure}
\begin{center}
\psfrag{a}{$H_1$}\psfrag{b}{$H_2$}
\psfrag{c}{$\ldots$}
\psfrag{A}{$G'$}\psfrag{B}{$G''$}
\includegraphics[width=80mm]{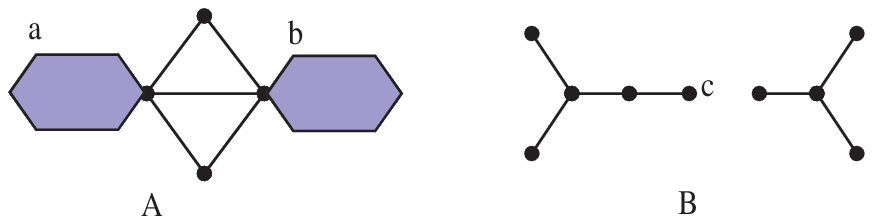}
  \caption{Graphs $G'$ and $G''$. }\label{fig1}
\end{center}
\end{figure}
\begin{lem}[\cite{Z-L-Z}]\label{lem:2.1}
Let $G$ be an $n$-vertex connected graph of girth $r=3$ and at least one block is non-complete. Then $Sz(G)-W(G)\geqslant 2,$ the equality holds if and only if $G\cong G',$ where $G'$ is depicted in Fig. \ref{fig1} satisfying each block of $H_1$ (resp. $H_2$) being a complete graph.
\end{lem}

\cite{T-D} showed that if $G$ is an $n$-vertex unicyclic graph with $n\geqslant 6,$ then
$W(G)\leqslant W(L_{n,3}).$ In fact, if $n=4,5,$ it is routine to check that this result also holds. Hence,
\begin{lem}\label{lem:2.2}
Let $G$ be a unicyclic graph on $n\geqslant 4$ vertices, then $W(G)\leqslant W(L_{n,3}).$
\end{lem}

For the unicyclic graph $C_r(T_1,T_2,\ldots,T_r),$ if the non-trivial rooted trees are just $T_i, T_j,\ldots, T_k$ with $i,j,\ldots, k\in \{1,2,\ldots, r\}$, then we denote it by $C_r^{i,j,\ldots, k}(T_i, T_j,\ldots, T_k).$ In particular, if every non-trivial tree $T_l,\, l=i,j,\ldots,k,$ is a star whose root is just its maximal degree vertex, then we denote it by $C_r^{i,j,\ldots, k}(S_{n_i+1}, S_{n_j+1},\ldots, S_{n_k+1});$ if every non-trivial tree $T_l, l=i,j,\ldots,k,$ is a path whose root is just an end-vertex, then we denote it by $C_r^{i,j,\ldots, k}(P_{n_i+1}, P_{n_j+1},\ldots, P_{n_k+1}).$
\begin{lem}[\cite{T-D}]\label{lem:2.7}
Let $\hat{G}=C_r^{i,j,\ldots, k}(S_{n_i+1},S_{n_j+1},\ldots,S_{n_k+1})$ and $\tilde{G}=\linebreak C_r^{i,j,\ldots, k}(P_{n_i+1},P_{n_j+1},\ldots,P_{n_k+1}).$ Then
$$
W(\hat{G})\leqslant W(G)\leqslant W(\tilde{G})
$$
for any graph $G=C_r^{i,j,\ldots, k}(T_i, T_j,\ldots, T_k)$ with $|T_t|=n_t+1, t=i,j,\ldots,k.$ The equality on the left (resp. right) holds if and only if $G\cong \hat{G}$ (resp. $G\cong \tilde{G}$).
\end{lem}
\begin{lem}[\cite{L-Z}]\label{lem:2.8}
Let $G$ be an $n$-vertex unicyclic graph with girth $r\geqslant 5$. If $r$ is odd, then
$
  W(G)\leqslant (n^3-25n+90)/6
$
with equality if and only if $G\cong L_{n,5}.$
\end{lem}

Let $T_n(n_1^{l_1},n_2^{l_2},\ldots,n_k^{l_k})$ be an $n$-vertex tree obtained by identifying one end-vertex for each of $l_1$ paths of length $n_1$,
$l_2$ paths of length $n_2,\ldots,$ $l_k$ paths of length $n_k.$ Clearly, $\sum_{i=1}^k l_in_i=n-1.$ For simplicity, let $A_{n,k}:=T_n(n-k-1, 1^k)$. A leaf in $A_{n,k}$ is called a \textit{unit leaf} if it is adjacent to the unique maximum degree vertex. Let $A_{n,k}^i$ be a graph obtained from $A_{n,k}$ by adding $i$ edges to connect its unit leaves.
\begin{lem}[\cite{D}]\label{lem:2.9}
Let $T$ be a tree of order $n\geqslant 9$. If $T\not\in\{ P_n,T_n(n-3,1^2),T_n(n-4,2,1),T_n(n-5,3,1),G''\},$ where $G''$ is depicted in Fig. \ref{fig1}. Then
$$
W(T)<W(T_n(n-5,3,1))<W(G'')<W(T_n(n-4,2,1))<W(T_n(n-3,1^2))<W(P_n).
$$
\end{lem}

By a direct calculation, we may obtain the Wiener index of each $n$-vertex trees for $5\leqslant n\leqslant 8$ (based on Table 2 of the Appendix in \cite{C-B-S}). Together with Lemma \ref{lem:2.9}, one obtains
\begin{lem}\label{lem:2.10}
Let $T$ be a tree of order $n\geqslant 5$. If $T\not\in\{P_n,T_n(n-3,1^2)\},$ then $W(T)<W(T_n(n-3,1^2))<W(P_n).$ Furthermore, if $n=8$ and $T\not\in \{P_8,T_8(5,1^2),T_8(4,2,1),T_8(3^2,1),G''\},$ then
$
W(T)<W(G'')<W(T_8(3^2,1))<W(T_8(4,2,1))<W(T_8(5,1^2))<W(P_n),
$
where $G''$ is the tree on $8$ vertices as depicted in Fig.~\ref{fig1}.
\end{lem}
\begin{lem}\label{lem:2.11}
Let $G$ be a graph obtained from vertex-disjoint connected graphs $G_1$ and $G_2$ by identifying a vertex of $G_1$ with a vertex of $G_2$ and denote the common vertex by $v.$ Then
$$
W(G)= W(G_1)+W(G_2)+(|G_2|-1)D_{G_1}(v)+(|G_1|-1)D_{G_2}(v).
$$
\end{lem}
\begin{proof}
Note that $d_G(x,y)=d_{G_i}(x,y)$ for $\{x,y\}\subseteq V_{G_i},i=1,2.$
By (\ref{eq:1.1}), we have
\begin{eqnarray}
W(G)&=&\sum_{x,y \in V_{G_1}}d_G(x, y)+\sum_{x,y \in V_{G_2}}d_G(x, y)
+\sum_{\substack{x\in V_{G_1}\backslash \{v\},\\ y\in V_{G_2}\backslash \{v\}}}d_G(x, y)\notag\\
&=&\sum_{x,y \in V_{G_1}}d_{G_1}(x, y)+\sum_{x,y \in V_{G_2}}d_{G_2}(x, y)
+\sum_{\substack{x\in V_{G_1}\backslash \{v\},\\ y\in V_{G_2}\backslash \{v\}}}(d_{G_1}(x,v)+d_{G_2}(v,y))\notag\\
&=& W(G_1)+W(G_2)+\sum_{\substack{x\in V_{G_1}\backslash \{v\},\\ y\in V_{G_2}\backslash \{v\}}}d_{G_1}(x,v)
+\sum_{\substack{x\in V_{G_1}\backslash \{v\},\\ y\in V_{G_2}\backslash \{v\}}}d_{G_2}(y,v)\notag\\
&=& W(G_1)+W(G_2)+(|G_2|-1)\sum_{x\in V_{G_1}\backslash \{v\}}d_{G_1}(x,v)
+(|G_1|-1)\sum_{y\in V_{G_2}\backslash \{v\}}d_{G_2}(y,v)\notag\\
&=&W(G_1)+W(G_2)+(|G_2|-1)D_{G_1}(v)+(|G_1|-1)D_{G_2}(v),\notag
\end{eqnarray}
as desired.
\end{proof}

The following result is a direct consequence of the above lemma.
\begin{cor}\label{cor2.1}
Let $u$ be a pendant vertex of an $n$-vertex connected graph $G$ and $v$ be the unique neighbor of $u$. Then
$
W(G)= W(G-u)+D_{G-u}(v)+n-1.
$
\end{cor}

\begin{lem}\label{lem:2.12}
Given a connected graph $H$ containing at least one edge, 
let $G_1$ be a graph obtained by identifying a vertex, say $v,$ of $H$ with a vertex in $C_k,$
$G_2$ be a graph obtained by identifying the vertex $v$ of $H$ with a minimal degree vertex in $L_{k,3}.$
Then one has $W(G_1)\leqslant W(G_2)$ with equality if and only if $k=3.$
\end{lem}
\begin{proof}
Clearly, if $k=3,$ then $G_1\cong G_2$. So we consider $k\geqslant 4$ in what follows. By Lemma \ref{lem:2.11}, we have
\begin{eqnarray}
W(G_1)&=& W(H)+W(C_k)+(k-1)D_{H}(v)+(|H|-1)D_{C_k}(v),\notag\\
W(G_2)&=& W(H)+W(L_{k,3})+(k-1)D_{H}(v)+(|H|-1)D_{L_{k,3}}(v).\notag
\end{eqnarray}
Thus,
\begin{eqnarray}
W(G_1)-W(G_2)=W(C_k)-W(L_{k,3})+(|H|-1)(D_{C_k}(v)-D_{L_{k,3}}(v)).\label{eq2.1}
\end{eqnarray}
It is easy to see that $|H|\geqslant 2.$ On the one hand, by Lemma \ref{lem:2.2}, we have $W(C_k)\leqslant W(L_{k,3}).$ On the other hand, $D_{C_k}(v)-D_{L_{k,3}}(v)\leqslant\frac{k^2}{4}-\frac{k^2-k-2}{2}=\frac{-k^2+2k+4}{4}<0$
for $k\geqslant 4.$ In view of (\ref{eq2.1}), we have $W(G_1)-W(G_2)<0,$ i.e., $W(G_1)<W(G_2),$  as desired.
\end{proof}

\begin{lem}\label{lem:2.13}
Let $G$ be an $n$-vertex connected graph of diameter at most $n-3.$ For all $v\in V_G,$ one has
$
D_G(v)\leqslant (n^2-n-6)/2
$
with equality if and only if $G\in \{A_{n,4}^i: i=0,1,2,3\}$ and $v$ is a non-unit leaf in $A_{n,4}^i, 0\leqslant i\leqslant 3$.
\end{lem}
\begin{proof}
Let $d$ be the diameter of $G,$ then by the definition of $D_G(v)$ for all $v\in V_G,$ we have
\begin{eqnarray}
D_G(v)&=& \sum_{x\in V_G}d_G(x,v)\notag\\
        &\leqslant& (1+2+\cdots+d)+d(n-d-1)\label{eq:2.5}\\
        &=& -\frac{1}{2}\left[\left(d-\frac{2n-1}{2}\right)^2-\frac{(2n-1)^2}{4}\right]\notag\\
        &\leqslant& \frac{n^2-n-6}{2}.\ \ \ \ \ \ \ \ \text{(Since $d\leqslant n-3$)}\label{eq:2.6}
\end{eqnarray}
The equality in (\ref{eq:2.5}) holds if and only if there are exactly $n-d$ vertices each of which is of distance $d$ from $v,$ whereas the equality in (\ref{eq:2.6}) holds if and only if $d=n-3.$ Hence, $D_G(v)=\frac{n^2-n-6}{2}$ if and only if $d=n-3, v$ is a non-unit leaf in $G$ and $G$ contains at most 3 unit leaves, i.e., $G\in \{A_{n,4}^i: i=0,1,2,3\}$.
\end{proof}
\begin{lem}\label{lem:2.14}
For the lollipop $L_{n,r}$ with $r<n$, if $r$ is even then $$W(L_{n,r})=\frac{1}{6}\left[n^3+\left(-\frac{3}{2}r^2+3r-1\right)n+
\left(\frac{5}{4}r^3-3r^2+r\right)\right].$$
\end{lem}
\begin{proof}
For convenience, let $u$ be the unique vertex of degree 3 in $L_{n,r}.$
By a direct calculation, we obtain that $W(C_r)=\frac{r^3}{8},$ $W(P_{n-r+1})=\frac{(n-r+1)(n-r)(n-r+2)}{6},$
$D_{C_r}(u)=\frac{2W(C_r)}{r}=\frac{r^2}{4}$ and $D_{P_{n-r+1}}(u)=1+2+\cdots+(n-r)=\frac{(n-r)(n-r+1)}{2}.$
Note that $u$ is a cut vertex of $L_{n,r}.$ Hence, we obtain (based on Lemma \ref{lem:2.11})
\begin{eqnarray*}
W(L_{n,r})&=& W(C_r)+W(P_{n-r+1})+(r-1)D_{P_{n-r+1}}(u)+(n-r)D_{C_r}(u)\\
        &=& \frac{r^3}{8}+\frac{(n-r+1)(n-r)(n-r+2)}{6}+\frac{(r-1)(n-r)(n-r+1)}{2}+\frac{r^2(n-r)}{4}\\
        &=& \frac{1}{6}\left[n^3+\left(-\frac{3}{2}r^2+3r-1\right)n+
\left(\frac{5}{4}r^3-3r^2+r\right)\right],
\end{eqnarray*}
as desired.
\end{proof}

Recall that $C_r(P_{n_1+1},P_{n_2+1},\ldots, P_{n_r+1})$ is an $n$-vertex unicyclic graph obtained from $C_r=v_1v_2\ldots v_rv_1$ and paths $P_{n_i+1}$'s by identifying an end-vertex of $P_{n_i+1}$ with $v_i$ on the cycle $C_r,\, i=1,2,\ldots, r.$
\begin{lem}\label{lem:3.1}
Let $G=C_r(P_{n_1+1},\ldots,P_{n_k+1},\ldots,P_{n_t+1},\ldots, P_{n_r+1})$ be an $n$-vertex unicyclic graph containing at least two non-trivial rooted paths, say $P_{n_k+1}$ and $P_{n_t+1}.$
\begin{wst}
\item[{\rm (i)}] If $\sum\limits_{\substack{i=1\\i\neq k,t}}^r(n_i+1)d_G(v_i,v_k)\leqslant\sum\limits_{\substack{i=1\\i\neq k,t}}^r(n_i+1)d_G(v_i,v_t),$
then
$$
W(G)<W(C_r(P_{n_1+1},\ldots,P_{n_{k-1}+1},P_1,P_{n_{k+1}+1}\ldots,
P_{n_{t-1}+1},P_{n_k+n_t+1},P_{n_{t+1}+1}\ldots,P_{n_r+1})).
$$
\item[{\rm (ii)}] If $\sum\limits_{\substack{i=1\\i\neq k,t}}^r(n_i+1)d_G(v_i,v_k)>\sum\limits_{\substack{i=1\\i\neq k,t}}^r(n_i+1)d_G(v_i,v_t),$
then
$$
W(G)<W(C_r(P_{n_1+1},\ldots,P_{n_{k-1}+1},P_{n_k+n_t+1},P_{n_{k+1}+1}\ldots,
P_{n_{t-1}+1},P_1,P_{n_{t+1}+1}\ldots,P_{n_r+1})).
$$
\end{wst}
\end{lem}
\begin{proof}
For convenience, denote by $P_{n_k+1}=u_1u_2\ldots u_{n_k+1}$ and $P_{n_t+1}=w_1w_2\ldots w_{n_t+1},$ where $u_1=v_k,w_1=v_t.$

(i) Let $G_1=G-u_1u_2+u_2w_{n_t+1},$ i.e.,
$$
  G_1=C_r(P_{n_1+1},\ldots,P_{n_{k-1}+1},P_1,P_{n_{k+1}+1}\ldots,
P_{n_{t-1}+1},P_{n_k+n_t+1},P_{n_{t+1}+1}\ldots,P_{n_r+1}).
$$
In what follows, we show that $W(G)<W(G_1).$

Let $A=(V_G\backslash V_{P_{n_k+1}})\bigcup \{u_1\}.$ Clearly, $A=(V_{G_1}\backslash V_{P_{n_k+1}})\bigcup \{u_1\}.$ By the definition of $W(G),$ we have
\begin{eqnarray}
W(G)&=&W(P_{n_k+1})+\sum_{x,y\in A} d_G(x,y)+\sum_{\substack{x\in V_{P_{n_k+1}}\backslash\{u_1\},\\y\in A\backslash\{u_1\}}} d_G(x,y)\notag\\
   &=&W(P_{n_k+1})+\sum_{x,y\in A} d_G(x,y)+\sum_{\substack{x\in V_{P_{n_k+1}}\backslash\{u_1\},\\ y\in A\backslash\{u_1\}}} (d_G(x,v_k)+d_G(v_k,y))\notag\\
&=&W(P_{n_k+1})+\sum_{x,y\in A} d_G(x,y)+(n-n_k-1)\sum_{x\in V_{P_{n_k+1}}\backslash\{u_1\}} d_G(x,v_k)+n_k\sum_{y\in A}d_G(v_k,y).\notag
\end{eqnarray}
Similarly, we have
\begin{align*}
W(G_1)=&\ W(P_{n_k+1})+\sum_{x,y\in A} d_{G_1}(x,y)+(n-n_k-1)\sum_{x\in V_{P_{n_k+1}}\backslash\{u_1\}} d_{G_1}(x,w_{n_t+1})\\
&+n_k\sum_{y\in A}d_{G_1}(w_{n_t+1},y).
\end{align*}
Note that $d_G(x,y)=d_{G_1}(x,y),$ $d_G(w_{n_t+1},y)=d_{G_1}(w_{n_t+1},y)$ and $d_G(z,v_k)=d_{G_1}(z,w_{n_t+1})$ for all $x,y\in A,  z\in V_{P_{n_k+1}}\backslash\{u_1\}.$ Hence,
\begin{eqnarray}\label{eq:3.01}
W(G)-W(G_1)=n_k\sum_{y\in A}(d_G(v_k,y)-d_G(w_{n_t+1},y)).
\end{eqnarray}
In what follows, we show that the right of (\ref{eq:3.01}) is negative. By a direct calculation, we have
\begin{eqnarray}
\sum_{y\in A}d_G(v_k,y)&=&\sum\limits_{\substack{i=1\\i\neq k}}^r\sum_{y\in V_{P_{n_i+1}}}(d_G(v_k,v_i)+d_G(v_i,y))\notag\\
&=&(n_t+1)d_G(v_k,v_t)+\sum\limits_{\substack{i=1\\i\neq k,t}}^r(n_i+1)d_G(v_k,v_i)+\sum\limits_{\substack{i=1\\i\neq k}}^rD_{P_{n_i+1}}(v_i)\notag
\end{eqnarray}
and
\begin{align}
\sum_{y\in A}d_G(w_{n_t+1},y)=&n_t+d_G(v_k,v_t)+\sum\limits_{\substack{i=1\\i\neq k,t}}^r\sum_{y\in V_{P_{n_i+1}}}(n_t+d_G(v_t,v_i)+d_G(v_i,y))+D_{P_{n_t+1}}(v_t)\notag\\
=&(n-n_k-n_t-1)n_t+d_G(v_k,v_t)+\sum\limits_{\substack{i=1\\i\neq k,t}}^r(n_i+1) d_G(v_t,v_i)+\sum\limits_{\substack{i=1\\i\neq k}}^rD_{P_{n_i+1}}(v_i).\notag
\end{align}
Bearing in mind the condition in (i) 
we have
\begin{equation}\label{eq:3.02}
\sum_{y\in A}d_G(v_k,y)-\sum_{y\in A}d_G(w_{n_t+1},y)\leqslant  n_t(d_G(v_k,v_t)-n+n_k+n_t+1).
\end{equation}

If $d_G(v_k,v_t)=1,$ then $d_G(v_k,v_t)-n+n_k+n_t+1=n_k+n_t+2-n=-\sum_{i\neq k,t}(n_i+1)\leqslant  -1.$ Together with (\ref{eq:3.01}) and (\ref{eq:3.02}), we have $W(G)<W(G_1)$.

If $d_G(v_k,v_t)\geqslant 2,$ then one may observe that $n-n_k-n_t-2=\sum_{i\neq k,t}(n_i+1)\geqslant 2d_G(v_k,v_t)-2.$ Thus, $d_G(v_k,v_t)-n+n_k+n_t+1\leqslant 1-d_G(v_k,v_t)\leqslant -1.$ Thus, in view of (\ref{eq:3.02}) we have $\sum_{y\in A}d_G(v_k,y)<\sum_{y\in A}d_G(w_{n_t+1},y).$ Together with (\ref{eq:3.01}), we obtain $W(G)<W(G_1),$ as desired.

(ii)\ By the same discussion as in the proof of (i), we may show that (ii) holds, which is omitted here.
\end{proof}

\begin{lem}\label{lem:3.2}
Let $C_r(T_1,T_2,\ldots,T_r)$ be an $n$-vertex unicyclic graph with girth $r\geqslant 6$. If $r$ is even, then
$$
  W(C_r(T_1,T_2,\ldots,T_r))\leqslant \frac{n^3-37n+168}{6}
$$
with equality if and only if $C_r(T_1, T_2,\ldots,T_r)\cong L_{n,6}$ if $n\not=8$ and $C_r(T_1,T_2,\ldots,T_r)\cong L_{8,6}$ or $L_{8,8},$
otherwise.
\end{lem}
\begin{proof}
Repeated using Lemmas \ref{lem:2.7} and \ref{lem:3.1} yields
\begin{align}
W(C_r(T_1,T_2,\ldots,T_r))\leqslant&\ W(L_{n,r})\label{eq:3.1}\\
                          =&\ \frac{1}{6}\left[n^3+\left(-\frac{3}{2}r^2+3r-1\right)n+
\left(\frac{5}{4}r^3-3r^2+r\right)\right]. \tag{By Lemma \ref{lem:2.14}}
\end{align}
The equality in (\ref{eq:3.1}) holds if and only if $G\cong L_{n,r}.$

In order to complete the proof, it suffices to compare the Wiener index of $L_{n,6}$ with that of $L_{n,r}$ for $r\geqslant 8.$ In fact, if $r\geqslant 8$ then by a direct calculation we have
\begin{align}
W(L_{n,6})-W(L_{n,r})&=\frac{1}{6}\left[\left(\frac{3}{2}r^2-3r-36\right)n-
\left(\frac{5}{4}r^3-3r^2+r-168\right)\right]\notag\\
&\geqslant \frac{1}{6}\left[\left(\frac{3}{2}r^2-3r-36\right)r-
\left(\frac{5}{4}r^3-3r^2+r-168\right)\right] \label{eq:3.5}\\
&=\frac{1}{6}\left(\frac{1}{4}r^3-37r+168\right) \notag\\
&\geqslant 0,\label{eq:3.6}
\end{align}
where (\ref{eq:3.5}) follows by the fact that $n\geqslant r$ and $\frac{3}{2}r^2-3r-36>0$ for $r\geqslant 8;$  the equality in (\ref{eq:3.5}) holds if and only if $n=r;$ whereas (\ref{eq:3.6}) follows by $r\geqslant 8;$ the equality in (\ref{eq:3.6}) holds if and only if $r=8.$ Hence, if $r\geqslant 8,$ then $W(L_{n,6})=W(L_{n,r})$ holds if and only if $n=r=8.$ That is to say, $W(L_{8,6})=W(L_{8,8}).$

By Lemma 3.12, one has $L_{n,6}=\frac{n^3-37n+168}{6}.$ Hence, together with (\ref{eq:3.1})--(\ref{eq:3.6}), it follows that
$W(C_r(T_1,\linebreak T_2, \ldots,T_r))\leqslant \frac{n^3-37n+168}{6}.$ The equality holds if and only if $C_r(T_1, T_2,\ldots,T_r)\cong L_{n,6}$ if $n\not=8,$ and $C_r(T_1,T_2,\ldots,T_r)\cong L_{8,6}$ or $L_{8,8}$ if $n=8,$ as desired.
\end{proof}
\begin{lem}\label{lem:3.5}
Let $G$ be a bicyclic graph of order $n$ with diameter $n-2.$ Then
$$
\frac{Sz(G)}{W(G)}\geqslant 1+\frac{12(n-3)}{n^3-13n+30}
$$
with equality if and only if $G\cong B(n,0).$
\end{lem}
\begin{proof}
Note that the diameter of $G$ is $n-2,$ by Lemma \ref{lem:2.6} we obtain $G\cong B(n,s),$ where $0 \leqslant s \leqslant \lfloor\frac{n-4}{2}\rfloor.$
By a direct calculation, one has
\begin{eqnarray}
W(B(n,s))&=& W(P_{n-1})+(1+2+\cdots+(s+1))+(1+2+\cdots+(n-s-3))+1\notag\\
&=&\frac{(n-2)(n-1)n}{6}+\frac{(s+1)(s+2)}{2}+\frac{(n-s-3)(n-s-2)}{2}+1\label{eq2.2}
\end{eqnarray}
and
\begin{align}
Sz(B(n,s))=&((n-1)+2(n-2)+\cdots+s(n-s))+2(s+2)(n-s-3)+2(s+1)(n-s+2)\notag\\
&+1+((n-1)+2(n-2)+(n-s-4)(s+4))\notag\\
=&\frac{ns(s+1)}{2}-\frac{s(s+1)(2s+1)}{6}+(4s+6)n-4(s+2)^2+1\notag\\
&+\frac{n(n-s-4)(n-s-3)}{2}-\frac{(n-s-4)(n-s-3)(2n-2s-7)}{6}\notag\\
=&\frac{(n-1)n(n+1)}{6}+ns-s^2-4s-1.\label{eq2.9}
\end{align}
Note that $0 \leqslant s \leqslant \lfloor\frac{n-4}{2}\rfloor-1.$ Hence, $n\geqslant 2s+6.$
Based on (\ref{eq2.2}),  one has
$W(B(n,s+1))-W(B(n,s))=2s+5-n<0,$ i.e., $W(B(n,s+1))<W(B(n,s)).$
By (\ref{eq2.9}) one has
$Sz(B(n,s+1))-Sz(B(n,s))=n-(2s+5)>0,$ i.e., $Sz(B(n,s+1))>Sz(B(n,s)).$
Thus, we have
\begin{align}
&W(B(n,\lfloor\frac{n-4}{2}\rfloor))<W(B(n,\lfloor\frac{n-4}{2}\rfloor-1))
<\cdots<W(B(n,1))<W(B(n,0)),\label{eq:2.10}\\
&Sz(B(n,0))<Sz(B(n,1))<\cdots<Sz(B(n,\lfloor\frac{n-4}{2}\rfloor-1))
<Sz(B(n,\lfloor\frac{n-4}{2}\rfloor)),\notag
\end{align}
which implies that
\begin{eqnarray*}
\frac{Sz(B(n,0))}{W(B(n,0))}<\frac{Sz(B(n,1))}{W(B(n,1))}<\cdots
<\frac{Sz(B(n,\lfloor\frac{n-4}{2}\rfloor-1))}{W(B(n,\lfloor\frac{n-4}{2}\rfloor-1))}
<\frac{Sz(B(n,\lfloor\frac{n-4}{2}\rfloor))}{W(B(n,\lfloor\frac{n-4}{2}\rfloor))}.
\end{eqnarray*}
Based on Eqs.(\ref{eq2.2})-(\ref{eq2.9}), one has $W(B(n,0))=\frac{n^3-13n+30}{6}$ and
$Sz(B(n,0))=\frac{n^3-n-6}{6}.$
Hence, $\frac{Sz(G)}{W(G)}\geqslant \frac{Sz(B(n,0))}{W(B(n,0))}=\frac{n^3-n-6}{n^3-13n+30}=1+\frac{12(n-3)}{n^3-13n+30}.$
The equality holds if and only if $G\cong B(n,0).$
\end{proof}

\section{Proof of Theorem \ref{thm5.1} }

Let $\mathscr{U}_n$ be the set of all unicyclic graph on $n\geqslant 10$ vertices. Recall that $\mathscr{A}_n=\{L_{n,3}, L_{n,4}, H_n^0,H_n^1,H_n^2, H_n^3, \linebreak C_3(P_{n-3}, P_2,P_1),C_3(P_{n-4},P_3,P_1),C_4(P_{n-4},P_1,P_2,P_1)\}.$ By a direct calculation, for $n\geqslant 10$, one has
\begin{align}
&W(H_n^3)=W(H_n^2)<W(H_n^1)<W(C_4(P_{n-4},P_1,P_2,P_1))<W(C_3(P_{n-4},P_3,P_1))<W(H_n^0)\notag\\
&\ \ \ \ <W(C_3(P_{n-3},P_2,P_1))=W(L_{n,4})<W(L_{n,3}).
\end{align}
In particular, $W(H_{10}^1)=W(C_3(P_5,P_4,P_1))$ and $W(H_{11}^2)=W(H_{11}^3)=W(C_3(P_6,P_4,P_1)).$

Note that $W(H_n^2)=W(H_n^3)=\frac{n^3-19n+54}{6}.$ Hence, in order to complete the proof, it suffices to show that $W(G)< \frac{n^3-19n+54}{6}$ for $G$ in $\mathscr{U}_n\setminus \mathscr{A}_n$ if $n\geqslant 12$ and for $G$ in $\mathscr{U}_n\setminus (\mathscr{A}_n\bigcup \{C_3(P_{n-5},P_4,P_1)\})$ if $n=10, 11.$

By Lemmas \ref{lem:2.8} and \ref{lem:3.2}, one has $W(G)<\frac{n^3-19n+54}{6}$ for all unicyclic graph $G$ of girth $r\geqslant 5.$ Hence, it suffices to consider $r=3,\,4.$
Choose such an $n$-vertex unicyclic graph $G$ of girth $r\leqslant 4$ in $\mathscr{U}_n\setminus \mathscr{A}_n$ if $n\geqslant 12$ and in $\mathscr{U}_n\setminus (\mathscr{A}_n\bigcup \{C_3(P_{n-5},P_4,P_1)\})$ if $n=10, 11$ such that $W(G)$ is as large as possible.

In what follows we show that $G$ does not contain at least three non-trivial rooted trees. Otherwise, assume without loss of generality that $T_1, T_k$ are the non-trivial rooted trees with the first two smallest sizes. Recall that for each rooted tree $T_i$ in $G$, one has $|T_i|=n_i+1$. Hence, by Lemma {\ref{lem:2.7}}, we have $G\cong C_r(P_{n_1+1},\ldots, P_{n_k+1},\ldots,P_{n_r+1}).$
By Lemma \ref{lem:3.1}, we have $W(G)<W(H')$ if $\sum_{ i\neq1,k}(n_i+1)d_G(v_i,v_1)\leqslant \sum_{i\neq1,k}(n_i+1)d_G(v_i,v_k)$ and $W(G)<W(H'')$ otherwise, where
\begin{eqnarray*}
H'&=&C_r(P_1,\ldots,P_{n_{k-1}+1},P_{n_1+n_k+1},P_{n_{k+1}+1},\ldots,P_{n_r+1}),\\
H''&=&C_r(P_{n_1+n_k+1},\ldots,P_{n_{k-1}+1},P_1,P_{n_{k+1}+1},\ldots,P_{n_r+1}).
\end{eqnarray*}
If both $H'$ and $H''$ are in $\mathscr{U}_n\setminus \mathscr{A}_n$ if $n\geqslant 12$ or in $\mathscr{U}_n\setminus (\mathscr{A}_n\bigcup \{C_3(P_{n-5},P_4,P_1)\})$ if $n=10, 11,$ then
 we obtain a contradiction to the maximality of $W(G).$
Otherwise, if $H'$ or $H''\in\mathscr{A}_n,$ as $H'$ and $H''$ have at least two non-trivial rooted trees, then $H'$ or $H''\in\{C_3(P_{n-3},P_2,P_1),C_3(P_{n-4},P_3,P_1),C_4(P_{n-4},P_1,\linebreak P_2,P_1)\}.$
Note that $T_1, T_k$ are the non-trivial rooted trees with the first two smallest sizes and $n_1+n_k\geqslant 2,$ we have $H'$ or $H''$ must be graph $C_3(P_{n-4},P_3,P_1).$ Together with the fact
$n\geqslant 10,$ we have $G\cong C_3(P_{n-4},P_2,P_2).$ By a simple computing, we have $W(G)=\frac{n^3-25n+96}{6}<\frac{n^3-19n+54}{6};$ if $H'$ or $H''\cong C_3(P_5,P_4,P_1),$ then
$G$ must be the graph $C_3(P_5,P_3,P_2), C_3(P_4,P_4,P_2)$ or $C_3(P_4,P_3,P_3).$
As $W(C_3(P_5,P_3,P_2))=135<144, W(C_3(P_4,P_4,P_2))=133<144$ and $W(C_3(P_4,P_3,P_3))=129<144,$
we obtain that $W(G)<144=(10^3-190+54)/6;$
if $H'$ or $H''\cong C_3(P_6,P_4,P_1),$ then
$G$ must be the graph $C_3(P_6,P_3,P_2)$ or $C_3(P_4,P_4,P_3).$
Since $W(C_3(P_6,P_3,P_2))=184<196$ and $W(C_3(P_4,P_4,P_3))=176<196,$
we have $W(G)<196=(11^3-209+54)/6,$ as desired.

Therefore, we obtain that $G$ contains at most two non-trivial rooted trees. Bearing in mind that $r=3$ or $r=4$. Hence, we proceed by considering the following two possible cases.

{\bf Case 1.} $r=3.$ In this case, it suffices to consider the following two subcases.

{\bf Subcase 1.1.} $G$ contains just one non-trivial rooted tree. Assume, without loss of generality, that $T_1$ is the non-trivial rooted tree. By Lemma \ref{lem:2.11}, we have
\begin{eqnarray}\label{eq3.5}
W(G)=W(T_1)+2D_{T_1}(v_1)+2n-3.
\end{eqnarray}
Note that $G\not\in\{ L_{n,3}, H_n^0,H_n^1\}.$ Hence, we have $T_1\not\in\{ P_{n-2},T_{n-2}(n-5,1^2),T_{n-2}(n-6,2,1)\}.$ This implies that the diameter of $T_1$ is at most $n-4.$

If $d(T_1)=n-4,$ as $T_1\not\in\{T_{n-2}(n-5,1^2),T_{n-2}(n-6,2,1)\},$ then $\eta_{T_1}^{n-4}(v_1)\leqslant 1,\,\eta_{T_1}^{n-5}(v_1)\leqslant 1.$ Thus,
\begin{eqnarray*}
D_{T_1}(v_1)&=&\eta_{T_1}^1(v_1)+2\eta_{T_1}^2(v_1)+\cdots+(n-5)\eta_{T_1}^{n-5}(v_1)+(n-4)\eta_{T_1}^{n-4}(v_1)\\
            &\leqslant &[1+2+\cdots+(n-5)+(n-4)]+(n-6)\left(\sum_{i=1}^{n-4}\eta_{T_1}^i(v_1)-n+4\right)\\
            &\leqslant &[1+2+\cdots+(n-5)+(n-4)]+(n-6)\ \ \ \ \ \ \ \ \ \ \ \text{(As $\sum_{i=1}^{n-4}\eta_{T_1}^i(v_1)=n-3$)}\\
            &\leqslant &\frac{n^2-5n}{2}.
\end{eqnarray*}
If $d(T_1)\leqslant n-5,$ by Lemma \ref{lem:2.13}, we have $D_{T_1}(v_1)\leqslant \frac{(n-2)^2-(n-2)-6}{2}=\frac{n^2-5n}{2}.$

Therefore, we obtain
\[\label{eq:4.0003}
  2D_{T_1}(v_1)\leqslant n^2-5n.
\]

If $|G|=10,$ then $|T_1|=8.$ Note that $T_1\not\in\{P_8,T_8(5,1^2),T_8(4,2,1)\}.$ By Lemma \ref{lem:2.10}, we have $W(T_1)\leqslant W(T_8(3^2,1))=75.$
By (\ref{eq3.5}) and (\ref{eq:4.0003}),  one has $W(G)\leqslant 75+50+17=142<144.$ Hence, (i) holds in this subcase.

Now we consider that $|G|\geqslant 11.$ If $T_1\cong G''$ ($G''$ is depicted in Fig. \ref{fig1}), then
$2D_{T_1}(v_1)\leqslant (n-5)(n-4)+2(n-3)=n^2-7n+14.$ By a direct calculation (based on (\ref{eq3.5})), we have $W(G)\leqslant\frac{n^3-31n+120}{6}<\frac{n^3-19n+54}{6}$ for $n\geqslant 11.$
If $T_1\ncong G'',$ by Lemma \ref{lem:2.9} we have $W(T_1)\leqslant W(T_{n-2}(n-7,3,1))=\frac{n^3-6n^2-7n+120}{6}.$
Thus, by (\ref{eq3.5}), we have
\begin{eqnarray*}
W(G)&\leqslant & \frac{n^3-6n^2-7n+120}{6}+(n^2-5n)+2n-3\\
&=&\frac{n^3-25n+102}{6}\\
&<&\frac{n^3-19n+54}{6}.\ \ \ \ \ \ \ \text{(As $n\geqslant 11$)}
\end{eqnarray*}
Hence, (ii) and (iii) hold in this subcase.

{\bf Subcase 1.2.} $G$ contains just two non-trivial rooted trees, say $T_1$ and $T_2.$ Assume, without loss of generality, that $|T_1|\geqslant |T_2|.$

If $n_2=1,$ then $n_1=n-4.$ As $G\ncong C_3(P_{n-3},P_2,P_1),$ we obtain that $T_1\ncong P_{n-3}.$
By Lemma \ref{lem:2.10}, we have $W(T_1)\leqslant W(T_{n-3}(n-6,1^2))=\frac{(n-5)(n-4)}{6}+2$ with equality if and only if $T_1\cong T_{n-3}(n-6,1^2).$ Observe that the diameter of $T_1$ is at most $n-5.$ Hence,
\begin{eqnarray}
D_{T_1}(v_1)&=& \eta_{T_1}^1(v_1)+2\eta_{T_1}^2(v_1)+\cdots+(n-5)\eta_{T_1}^{{n-5}}(v_1)\notag\\
&\leqslant &[1+2+3+\cdots+(n-5)]+\left(\sum_{i=1}^{n-5}\eta_{T_1}^i(v_1)-n+5\right)(n-5)\label{eq:1}\\
&=&[1+2+3+\cdots+(n-5)]+(n-5) \ \ \ \ \ \ \ \ \ \ \ \ \ \ \ \ \ \text{(As $\sum_{i=1}^{n-5}\eta_{T_1}^i(v_1)=n-4$)}\notag\\
&=&\frac{(n-5)(n-2)}{2},\label{eq:2}
\end{eqnarray}
 where the equality in (\ref{eq:1}) holds if and only if $\eta_{T_1}^1(v_1)=\cdots=\eta_{T_1}^{n-6}(v_1)=1$ and $\eta_{T_1}^{n-5}(v_1)=2.$ This means that $T_1\cong T_{n-3}(n-6,1^2)$ and $v_1$ is a non-unit leaf in $T_{n-3}(n-6,1^2).$ As $v_1$ is a cut vertex of $G,$ we obtain (based on Lemma \ref{lem:2.11})
\begin{eqnarray}
W(G)&=&W(T_1)+W(L_{4,3})+3D_{T_1}(v_1)+(n-4)D_{L_{4,3}}(v_1)\notag\\
    &\leqslant &\frac{(n-5)(n-4)}{6}+2+8+3D_{T_1}(v_1)+4(n-4)\label{eq:13.1}\\
    &\leqslant &\frac{(n-5)(n-4)}{6}+2+8+\frac{3(n-5)(n-2)}{2}+4(n-4)\label{eq:13.2}\\
    &=&\frac{n^3-19n+54}{6},\notag
\end{eqnarray}
where the equality in (\ref{eq:13.1}) holds if and only if $T_1\cong T_{n-3}(n-6,1^2);$ the equality in (\ref{eq:13.2}) holds if and only if $T_1\cong T_{n-3}(n-6,1^2)$ and $v_1$ is a non-unit leaf in $T_{n-3}(n-6,1^2).$ Thus, $G\cong H_n^3.$ As $G\not\in \mathscr{A}_n,$ we have $G\ncong H_n^3,$ which implies that $W(G)<\frac{n^3-19n+54}{6}.$

If $n_2=2,$ then $n_1=n-5$ and $T_2\cong P_3.$ This implies that $d_G(v_2)=3$ or $4.$ If $d_G(v_2)=3,$ by a similar discussion as in the proof of $n_2=1,$ we can obtain that $W(G)\leqslant\frac{n^3-25n+90}{6}
<\frac{n^3-19n+54}{6}$ for $n\geqslant 10.$ So we only consider the case $d_G(v_2)=4.$
Note that $|T_1|=n-4.$ By Lemma \ref{lem:2.10}, we have $W(T_1)\leqslant W(P_{n-4})=\frac{(n-4)(n-5)(n-3))}{6}.$ It is easy to see that $D_{T_1}(v_1)\leqslant 1+2+\cdots+(n-5)=\frac{(n-5)(n-4)}{2}.$ As $v_1$ is a cut vertex of $G,$ by Lemma \ref{lem:2.11} we have
\begin{eqnarray*}
W(G)&=&W(T_1)+W(G[V_{T_2}\cup \{v_1,v_3\}])+4D_{T_1}(v_1)+(n-5)D_{G[V_{T_2}\cup \{v_1,v_3\}]}(v_1)\notag\\
    &\leqslant &\frac{(n-4)(n-5)(n-3)}{6}+15+\frac{(n-5)(n-4)}{2}+6(n-5)\notag\\
    &=&\frac{n^3-25n+90}{6}\\
    &<&\frac{n^3-19n+54}{6}.\ \ \ \ \ \ \ \text{(Since $n\geqslant 10$)}\
\end{eqnarray*}

Now, we consider $n_2\geqslant 3.$  For $n=10,$ we obtain that $n_1=4, n_2=3$ directly. As $G\ncong C_3(P_5,P_4,P_1),$ we have $T_1\ncong P_5$ or $T_2\ncong P_4.$ If $T_1\ncong P_5,$ then by Lemma \ref{lem:2.10}, we have $W(T_1)\leqslant W(T_5(2,1^2))=18.$ Since $G[V_{T_2}\cup \{v_1,v_3\}]$ is a unicyclic graph on 6 vertices, by Lemma \ref{lem:2.2}, we have $W(G[V_{T_2}\cup \{v_1,v_3\}])\leqslant L_{6,3}=31.$ Observe that the diameters of $T_1$ and $T_2$ are at most $3,$ we can obtain that $D_{T_1}(v_1)\leqslant 1+2+3+3=9$ and $D_{G[V_{T_2}\cup \{v_1,v_3\}]}(v_1)\leqslant 1+1+2+3+4=11.$
By Lemma \ref{lem:2.11}, we have
\begin{align*}
\ \ \ W(G)=&\ W(T_1)+W(G[V_{T_2}\cup \{v_1,v_3\}])+3D_{T_1}(v_1)+4D_{G[V_{T_2}\cup \{v_1,v_3\}]}(v_1)\\
          \leqslant&\ 18+31+27+44=120<144.\notag
\end{align*}
If $T_1\cong P_5,$ then one must have $T_2\ncong P_4.$ Thus, $T_2\cong S_4$ and $d_G(v_2)=5$ or $d_G(v_2)=3.$ If $d_G(v_2)=5$,  then $W(G)=126<144;$ if $d_G(v_2)=5,$ then $W(G)=138<144.$ Hence, (i) holds in this subcase.

For $n\geqslant 11,$ by Lemma {\ref{lem:2.7}}, we have $G=C_3(P_{n_1+1},P_{n_2+1},P_1).$ Based on the structure of $G$ and the fact $n_1=n-3-n_2,$ we have
\begin{eqnarray}
W(G)&=&W(P_{n_1+n_2+2})+[1+2+\cdots+(n_1+1)]+[1+2+\cdots+(n_2+1)]\notag\\
    &=&W(P_{n-1})+[1+2+\cdots+(n-n_2-2)]+[1+2+\cdots+(n_2+1)]\notag \\
    &=&\frac{(n-1)(n-2)n}{6}+\frac{(n-n_2-2)(n-n_2-1)}{2}+\frac{(n_2+1)(n_2+2)}{2}\notag\\
    &=&n_2^2-(n-3)n_2+\frac{n^3-7n+12}{6}\notag\\
    &\leq& 9-3(n-3)+\frac{n^3-7n+12}{6} \ \ \ \ \ \ \text{(As $3\leqslant n_2\leqslant \left\lfloor\frac{n-3}{2}\right\rfloor$)}\label{eq3.71} \\
   &=& \frac{n^3-25n+120}{6}\notag \\
   &\leq& \frac{n^3-19n+54}{6}, \ \ \ \ \ \ \ \ \ \ \ \ \ \ \ \ \ \ \ \ \ \ \ \text{(As $n\geqslant 11$)}\label{eq3.8}
\end{eqnarray}
where the equality in (\ref{eq3.71}) holds if and only if $n_2=3$; whereas the equality in (\ref{eq3.8}) holds if and only if $n=11.$ That is, $G\cong C_3(P_6,P_4,P_1).$ As $G\ncong C_3(P_6,P_4,P_1)$ for $n=11,$ we have $W(G)<\frac{n^3-19n+54}{6}.$ Hence, (ii) and (iii) hold in this subcase.

Hence, by Subcases 1.1 and 1.2 we obtain that (i), (ii) and (iii) hold for Case 1.

{\bf Case 2.}\ $r=4.$ In this case, we first consider that $G$ contains just one non-trivial rooted tree, say $T_1$. By Lemma \ref{lem:2.11}, we have
\begin{eqnarray}\label{eq3.9}
W(G)=W(T_1)+3D_{T_1}(v_1)+4n-8.
\end{eqnarray}
 Note that $G\ncong L_{n,4},$ we have
$T_1\ncong P_{n-3}.$ Bearing in mind that $n\geqslant 10,$ together with Lemma \ref{lem:2.10} for $n-3\, (\geqslant 7),$ it follows that $W(T_1)\leqslant W(T_{n-3}(n-6,1^2))=\frac{n^3-9n^2+20n+12}{6}.$
By (\ref{eq:2}), we have $D_{T_1}(v_1)\leqslant \frac{(n-5)(n-2)}{2}$ with equality if and only if $T_1\cong T_{n-3}(n-6,1^2)$ and $v_1$ is a non-unit leaf in $T_{n-3}(n-6,1^2).$
Thus, by (\ref{eq3.9}), we have
\begin{eqnarray}
W(G)&\leqslant & \frac{n^3-9n^2+20n+12}{6}+3D_{T_1}(v_1)+4n-8\label{eq3.10}\\
&\leqslant &\frac{n^3-9n^2+20n+12}{6}+\frac{3(n-5)(n-2)}{2}+4n-8\label{eq3.11}\\
&=&\frac{n^3-19n+54}{6}.\notag
\end{eqnarray}
By Lemma \ref{lem:2.10}, the equality in (\ref{eq3.10}) holds if and only if $T_1\cong T_{n-3}(n-6,1^2);$ whereas the equality in (\ref{eq3.11}) holds if and only if $T_1\cong T_{n-3}(n-6,1^2)$ and $v_1$ is a non-unit leaf in $T_{n-3}(n-6,1^2).$ Hence, $G\cong H_n^2.$
As $G\not\in \mathscr{A}_n,$ we have $G\ncong H_n^2,$ which implies that $W(G)<\frac{n^3-19n+54}{6}.$

Now we consider the remaining case that $G$ contains just two non-trivial rooted trees.
By Lemma {\ref{lem:2.7}}, up to isomorphism one has $G=C_4(P_{n_1+1},P_{n_2+1},P_1,P_1)$ with $0\leqslant n_1\leqslant n_2,$ or $G=C_4(P_{n_1+1},P_1,\linebreak P_{n_3+1},P_1)$ with $0\leqslant n_1\leqslant n_3.$ For the former case, by a direct calculation, one has
$W(G)=2n_1^2-2(n-4)n_1+\frac{n^3-13n+36}{6}\leqslant 2-2(n-4)+\frac{n^3-13n+36}{6}=\frac{n^3-25n+96}{6}< \frac{n^3-19n+54}{6}$ for $n\geqslant 10.$
For the latter case, bearing in mind that $n_1\geqslant 2,$ one has
$W(G)=n_1^2-(n-4)n_1+\frac{n^3-13n+36}{6}\leqslant 2-2(n-4)+\frac{n^3-13n+36}{6}=\frac{n^3-25n+96}{6}<\frac{n^3-19n+54}{6}$ for $n\geqslant 10.$
Here, (i), (ii) and (iii) are true for Case 2.

This completes the proof.
\qed

\section{Proof of Theorem \ref{thm5.2} }
In this section, we give the proof of Theorem \ref{thm5.2}.
\begin{figure}[h!]
\begin{center}
\psfrag{a}{$F_1$}\psfrag{b}{$F_2$}\psfrag{c}{$F_3$}
\psfrag{d}{$F_4$}\psfrag{e}{$F_5$}
\psfrag{f}{$F_6$}\psfrag{g}{$F_7$}\psfrag{h}{$F_8$}
\psfrag{i}{$F_9$}\psfrag{j}{$F_{10}$}\psfrag{x}{$F_{11}$}
\includegraphics[width=130mm]{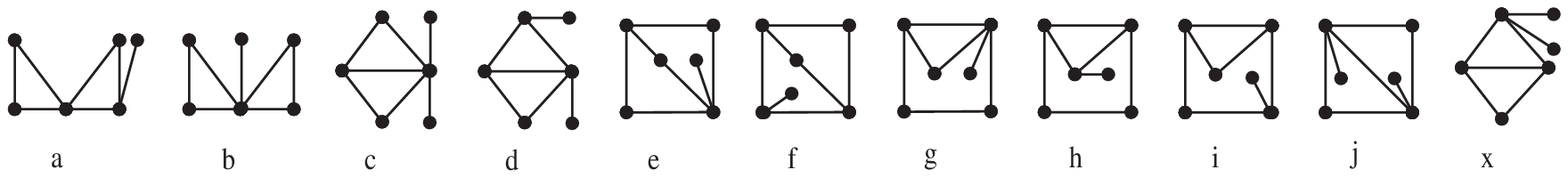}
  \caption{Graphs $F_i$ used in the proof of Theorem \ref{thm5.2},~$1\leqslant i\leqslant 11$.}\label{fig3}
\end{center}
\end{figure}
For convenience, let
$B_{k,l,t}$ be the bicyclic graph of order $n$ obtained from paths $P_{k+1},\, P_{l+1},\, P_{t+1}$ by identifying all the left (resp. right) end-vertices as a new vertex. One often calls $B_{k,l,t}$ the $\theta$-\textit{graph}.

\noindent{\bf Proof of Theorem \ref{thm5.2}.}\ \ Let $\mathscr{B}_n$ be the set of all bicyclic graph on $n\geqslant 6$ vertices. By a direct calculation, for $n\geqslant 7$, one has
\begin{align}
W(B_n^{(1)})<W(B(n,1))<W(B_n^{3,3})<W(B(n,0)).\notag
\end{align}
In particular, $W(B(8,2))=W(B_8^{(1)})$ and $W(B_6^{(1)})<W(B(6,1))=W(B_6^{3,3})<W(B(6,0)).$

It is routine to check that $W(B_n^{(1)})=(n^3-19n+54)/6.$ Hence, in order to complete the proof, it suffices to show that $W(G)< (n^3-19n+54)/6$ for all $G$ in $\mathscr{B}_n\setminus \{B(n,0),B(n,1),B_n^{(1)},B_n^{3,3}\}$ if $n\not=8,$ and for $G$ in $\mathscr{B}_8\setminus  \{B(8,0),B(8,1),B_8^{(1)},B_8^{3,3},B(8,2)\})$ if $n=8.$

Note that the diameter of $G$ is at most $n-2.$ If $d(G)=n-2,$ then $G\cong B(n,s)$ for some $0\leqslant s\leqslant \lfloor\frac{n-4}{2}\rfloor.$
As $G\not\in \{B(n,0),B(n,1)\},$ we have $s\geqslant 2,$ which implies that $n\geqslant 8.$ Note that $G\ncong B(8,2)$ for $n=8.$ Thus, $n\geqslant 9.$ Combining with (\ref{eq2.2}) and (\ref{eq:2.10}), we obtain $W(G)\leqslant W(B(n,2))=\frac{n^3-25n+102}{6}<\frac{n^3-19n+54}{6}$ for $n\geqslant 9.$
Hence, (ii) and (iii) hold for the case $d(G)=n-2.$

If $d(G)\leqslant n-3,$ then we show our result by induction on $n.$ When $n=6,$ we have $W(B_6^{(1)})=26$ and with the help of Nauty based on \cite{M-V0}, we obtain $\mathscr{B}_6\backslash \{B(6,0),B(6,1),B_6^{(1)},B_6^{3,3}\}\linebreak =\{B_6^{3,4},B_{1,2,4}, B_{1,3,3}, B_{2,2,3},
F_1,F_2,\ldots,F_{11}\},$ where $F_1, F_2,\ldots,F_{11}$ are depicted in Fig. \ref{fig3}. By a direct calculation, one has
$$
\begin{array}{lllll}
  W(F_1)=25,&W(F_4)=25,&W(F_7)=24,&W(F_{10})=24,
   &W(B_{1,2,4})=24,
\\[5pt]
W(F_2)=23,
&W(F_5)=24, &W(F_8)=25,&W(F_{11})=25,&W(B_{1,3,3})=25,
\\[5pt]
W(F_3)=23,&W(F_6)=25,
&W(F_9)=25,&W(B_6^{3,4})=25,&W(B_{2,2,3})=23.
\end{array}
$$
Thus, $W(G)\leqslant 25<26=W(B_6^{(1)}).$ That is, (i) is true for $n=6.$

 Now, consider $n\geqslant 7.$ Assume that our result is true for $n-1.$
 In order to complete the proof,
it suffices to consider the following two possible cases.
\vspace{2mm}

{\bf Case 1}.\ $G$ contains pendant vertices. Choose a pendant vertex, say $u$, in $G$ and let $v$ be its unique neighbor. Clearly, $G-u$ is a bicyclic graph on $n-1$ vertices. By Corollary \ref{cor2.1}, we have
\begin{eqnarray}
W(G)&=& W(G-u)+D_{G-u}(v)+n-1\label{eq3.7}
\end{eqnarray}
Note that the diameter $d(G)\leqslant n-3.$ Hence, $d(G-u)\leqslant n-3.$

If $d(G-u)=n-3,$ then by Lemma \ref{lem:2.6}, we have $G-u\cong B(n-1,s)$ with $0\leqslant s\leqslant \lfloor\frac{n-5}{2}\rfloor.$
Combining with (\ref{eq:2.10}), we obtain $W(G-u)\leqslant W(B(n-1,0))=\frac{n^3-3n^2-10n+42}{6}.$
Observe that $v$ isn't a pendant vertex of $G-u.$ Otherwise, $d(G)=n-2,$ a contradiction. Hence,
$D_{G-u}(v)\leqslant 1+2+3+\cdots+(n-4)+1+(n-5)=\frac{n^2-5n+4}{2}.$ By (\ref{eq3.7}), we have
\begin{eqnarray}
W(G)&\leqslant & \frac{n^3-3n^2-10n+42}{6}+\frac{n^2-5n+4}{2}+n-1\notag\\
        &=& \frac{n^3-19n+48}{6}\notag\\
        &<& \frac{n^3-19n+54}{6}.\notag
\end{eqnarray}

Now, consider that $d(G-u)\leqslant n-4.$ If $G-u\cong B_{n-1}^{3,3},$ then by a direct calculation, we have $W(G-u)=\frac{n^3-3n^2-10n+36}{6}.$ Based on the structure of $B_{n-1}^{3,3},$ we obtain $D_{G-u}(v)\leqslant 1+2+\cdots+(n-4)+1+(n-4)=\frac{n^2-5n+6}{2}.$ By (\ref{eq3.7}), one has
$W(G)\leqslant \frac{n^3-3n^2-10n+36}{6}+\frac{n^2-5n+6}{2}+n-1= \frac{n^3-19n+48}{6}<\frac{n^3-19n+54}{6},$ the result holds.
Note that $d(G-u)\leqslant n-4,$ by Lemma $\ref{lem:2.13},$ we have
\begin{eqnarray}
D_{G-u}(v)\leqslant \frac{(n-1)^2-(n-1)-6}{2}=\frac{n^2-3n-4}{2},\label{eq3.0}
\end{eqnarray}
the equality in (\ref{eq3.0}) if and only if $G-u\cong B_{n-1}^{(1)}$ and $v$ is a non-unit leaf in $B_{n-1}^{(1)}.$

If $G-u\cong B_{n-1}^{(1)},$ then $W(G-u)=\frac{n^3-3n^2-16n+72}{6}.$
By (\ref{eq3.7}) and (\ref{eq3.0}), one has
\begin{eqnarray}
W(G)&\leqslant & \frac{n^3-3n^2-16n+72}{6}+\frac{n^2-3n-4}{2}+n-1\label{eq2.4}\\
    &=& \frac{n^3-19n+54}{6}.\notag
\end{eqnarray}
By (\ref{eq3.0}), the equality in (\ref{eq2.4}) holds if and only if $G-u\cong B_{n-1}^{(1)}$ and $v$ is a non-unit leaf in  $B_{n-1}^{(1)}$. That is to say, $W(G)=\frac{n^3-19n+54}{6}$ holds if and only if $G\cong B_n^{(1)},$ which is impossible. Hence, $W(G)<\frac{n^3-19n+54}{6}.$

Thus, it is left for us to consider $G-u\not\in\{B_{n-1}^{3,3}, B_{n-1}^{(1)}\}.$ By the induction hypothesis, we have
$W(G-u)<\frac{(n-1)^3-19(n-1)+54}{6}=\frac{n^3-3n^2-16n+72}{6}.$
Combining with (\ref{eq3.7}) and (\ref{eq3.0}), we obtain $W(G)<\frac{n^3-3n^2-16n+72}{6}+\frac{n^2-3n-4}{2}+n-1=\frac{n^3-19n+54}{6},$ as desired. Therefore, (ii) and (iii) hold for Case 1.
 \vspace{2mm}

{\bf Case 2.}\ $G$ does not contain pendant vertices. In this case, $G$ is either the graph $B_n^{p,q}$ or the $\theta$-graph $B_{k,l,t}.$
So we proceed by distinguishing the following two possible subcases.
\vspace{2mm}

{\bf Subcase 2.1.}\ $G=B_n^{p,q}.$ Assume, without loss of generality, that $p\leqslant q.$ Note that $G\ncong B_n^{3,3},$ we have $q\geqslant 4.$
As $p\geqslant 3,$ by Lemma \ref{lem:2.12}, we have $W(G)\leqslant W(B_n^{3,q}).$
In what follows, we show that $W(B_n^{3,q})<\frac{n^3-19n+54}{6}.$

In fact, by Lemma \ref{lem:2.11}, one has
\begin{eqnarray}\label{eq2.5}
W(B_n^{3,q})=\left\{
  \begin{array}{ll}
    W(L_{n-2,q})+n^2-3n+3+\frac{-q^2+2q}{2}, & \hbox{ if $q$ is even;} \\
    W(L_{n-2,q})+n^2-3n+3+\frac{-q^2+2q-1}{2}, & \hbox{ if $q$ is odd}
  \end{array}
\right.
\end{eqnarray}
and using Lemma \ref{lem:2.14} for $n-2$ and $r=4,$ one has
\begin{eqnarray}\label{eq2.6}
W(L_{n-2,4})=\frac{(n-2)^2-13(n-2)+36}{6}=\frac{n^3-6n^2-n+54}{6}.
\end{eqnarray}
If $n=7,$ then $W(L_{5,q})\leqslant\max\{W(L_{5,4}),W(C_5)\}=16.$
By (\ref{eq2.5}), we have $W(B_7^{3,q})\leqslant 16+7^2-21+3-4=43<44=\frac{n^3-19n+54}{6},$ as desired. Now consider $n\geqslant 8.$
Note that $W(L_{n-2,5})=\frac{n^3-6n^2-13n+132}{6}<\frac{n^3-6n^2-n+54}{6}$ and $W(L_{n-2,6})=\frac{n^3-6n^2-25n+234}{6}<\frac{n^3-6n^2-n+54}{6}.$
Combining with (\ref{eq2.5})-(\ref{eq2.6}), we obtain
\begin{align}
W(B_n^{3,q})\leqslant & W(L_{n-2,q})+n^2-3n+3-4\tag{Since $q\geqslant 4$}\\
\leqslant & W(L_{n-2,4})+n^2-3n-1\tag{By Lemmas \ref{lem:2.8} and \ref{lem:3.2}}\\
=& \frac{n^3-19n+48}{6} \tag{By Eq.(\ref{eq2.6})}\\
<& \frac{n^3-19n+54}{6}.\notag
\end{align}
Hence, (ii) and (iii) hold for Subcase 2.1.
\vspace{2mm}

{\bf Subcase 2.2.}\ $G=B_{k,l,t}.$ In this subcase, $k+l+t=n+1.$ Assume, without loss of generality, that $k\leqslant l\leqslant t$ and $G$ contains just two vertices, say $u_1,u_2$, of degree 3. Clearly, each of the rest vertices is of degree 2.

Recall that $\eta_G^i(x)=|\{u | d_G(u,x)=i\}|$ for all $x\in V_G.$ As $d(G)\leqslant \lfloor\frac{n}{2}\rfloor,$ we have $\eta_G^i(x)=0$ for $i>\lfloor\frac{n}{2}\rfloor.$ In what follows, we show that if $\eta_G^{i+1}(x)> 0,$ then $\eta_G^i(x)\geqslant 2.$

In fact, if $\eta_G^i(x)<2,$ then $\eta_G^{i}(x)=1.$ This implies that there exists a unique vertex $u$ such that $d(u,x)=i.$ Since $\eta_G^{i+1}(x)> 0,$ we obtain that $u$ is a cut vertex of $G,$ a contradiction to the fact that $B_{k,l,t}$ contains no cut vertex. Hence, $\eta_G^i(x)\geqslant 2.$

Now we show that our result holds for even $n.$ Similarly, we can also show that our result holds for odd $n,$ which is omitted here.

Observe that, for all $x\in V_G\backslash \{u_1,u_2\},$ we have $\eta_G^{\frac{n}{2}}(x)\leqslant 1.$ Otherwise, $G$ has at least $n+1$ vertices, a contradiction. Similarly, we have $\eta_G^{\frac{n}{2}}(u_1)=\eta_G^{\frac{n}{2}}(u_2)=0,$ $\eta_G^{\frac{n-2}{2}}(u_1)\leqslant 2$ and $\eta_G^{\frac{n-2}{2}}(u_2)\leqslant 2.$
Thus, for all $x\in V_G\backslash \{u_1,u_2\},$ one has
\begin{eqnarray}\label{eq2.7}
D_G(x)=\sum_{u\in V_G}d_G(u,x)=\sum_{1\leqslant i\leqslant \frac{n}{2}} i\eta_G^i(x)\leqslant 2\left(1+2+\cdots+\frac{n-2}{2}\right)+\frac{n}{2}=\frac{n^2}{4}.
\end{eqnarray}
For $y\in \{u_1,u_2\},$ one has
\begin{eqnarray}\label{eq2.8}
D_G(y)=\sum_{u\in V_G}d_G(u,y)=\sum_{1\leqslant i\leqslant \frac{n-2}{2}} i\eta_G^i(y)\leqslant 3+2\left(2+\cdots+\frac{n-2}{2}\right)=\frac{n^2-2n+4}{4}.
\end{eqnarray}
Together with (\ref{eq:1.1}) and (\ref{eq2.7})-(\ref{eq2.8}), it follows that
\begin{eqnarray*}
W(G)&=&\frac{1}{2}\sum_{x\neq u_1,u_2}D_G(x)+\frac{1}{2}D_G(u_1)+\frac{1}{2}D_G(u_2)\\
&\leqslant & \frac{n-2}{2}\cdot \frac{n^2}{4}+\frac{n^2-2n+4}{4}\\
&=&\frac{n^3-4n+8}{8}\\
&<&\frac{n^3-19n+54}{6}.\ \ \ \ \ \ \ \ \ \text{(As $n\geqslant 7$)}
\end{eqnarray*}
Hence, (ii) and (iii) hold for Subcase 2.2.

This completes the proof.
\qed

\section{Proofs of Theorems \ref{thm4.1} and \ref{thm4.2}}
In this section,  we give the proofs of Theorems \ref{thm4.1} and \ref{thm4.2}.
\begin{figure}[h!]
\begin{center}
\psfrag{d}{$H_1$}\psfrag{e}{$H_2$}\psfrag{f}{$H_3$}
\psfrag{g}{$H_4$}\psfrag{h}{$H_5$}
\psfrag{j}{$H_6$}\psfrag{k}{$H_7$}\psfrag{l}{$H_8$}
\psfrag{m}{$H_9$}
\includegraphics[width=120mm]{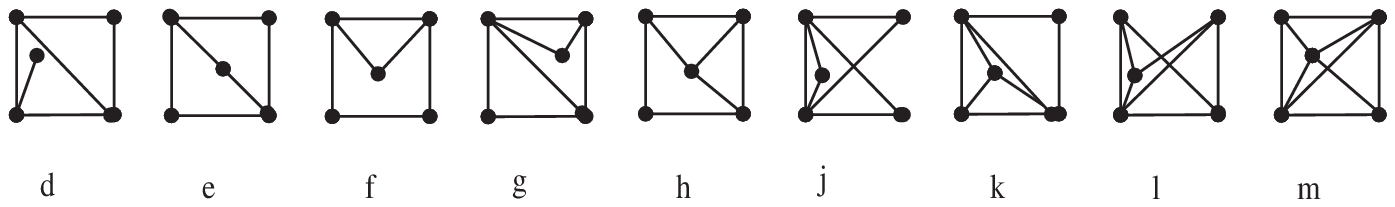}
  \caption{Graphs $H_i$ used in the proof of Theorem \ref{thm4.1},~$1\leqslant i\leqslant 9$.}
\end{center}
\end{figure}
\vspace{3mm}

\noindent{\bf Proof of Theorem \ref{thm4.1}.}\ \
If $n=5,$ then with the help of Nauty based on \cite{M-V0} all of such graphs are $C_5, L_{5,4},B_5^{(1)},H_1,H_2, \ldots, H_9,$
where $H_1, H_2, \ldots, H_9$ are depicted in Fig. 4. By a direct calculation, we have
$$
\begin{array}{lllllll}
 \frac{Sz(B_5^{(1)})}{W(B_5^{(1)})}=\frac{8}{7},
&\frac{Sz(C_5)}{W(C_5)}=\frac{4}{3},
&\frac{Sz(L_{5,4})}{W(L_{5,4})}=\frac{7}{4},
&\frac{Sz(H_1)}{W(H_1)}=\frac{19}{15},
&\frac{Sz(H_2)}{W(H_2)}=\frac{18}{7},
&\frac{Sz(H_3)}{W(H_3)}=\frac{12}{7},
\\[8pt]
\frac{Sz(H_4)}{W(H_4)}=\frac{18}{13},
&\frac{Sz(H_5)}{W(H_5)}=\frac{29}{13},
&\frac{Sz(H_6)}{W(H_6)}=\frac{19}{13},
&\frac{Sz(H_7)}{W(H_7)}=\frac{4}{3},
&\frac{Sz(H_{8})}{W(H_{8})}=2,
&\frac{Sz(H_{9})}{W(H_{9})}=\frac{15}{11}.
\\[5pt]
\end{array}
$$
Thus, $B_5^{(1)}$ is the unique extremal graph having minimum value on $\frac{Sz(G)}{W(G)}.$ Our result holds for $n=5.$

Now, we assume that $n\geqslant 6.$  Note that at least one block of $G$ is non-complete. Hence,
$G$ contains at least one cycle and its diameter $d(G)\leqslant n-2.$ For convenience, let $r$ be the girth of graph $G.$
Then if $r\geqslant 4,$ by Theorem \ref{lem:2.3}, we have
\begin{eqnarray}\label{eq3.1}
\frac{Sz(G)}{W(G)}&\geqslant& \min\left\{1+\frac{24(n-2)}{n^3-13n+36},1+\frac{6(2n-5)}{n^3-25n+90}\right\}\notag\\
&=&1+\frac{6(2n-5)}{n^3-25n+90}>1+\frac{12}{n^3-19n+54}.
\end{eqnarray}

So, in what follows, it suffices to consider graph $G$ of girth $r=3.$ In order to complete our proof, we distinguish the following
two cases $d(G)=n-2$ and $d(G)\leqslant n-3.$

If $d(G)=n-2,$ then there exists an induced path $P_{n-1}$ in $G.$ Note that $r=3$ and at least one block is non-complete, we obtain that $G$ contains only one vertex in $V_G\setminus V_{P_{n-1}}$ such that it must be adjacent to just three vertices of $P_{n-1}.$ Thus, $G$ is a bicyclic graph. By Lemma \ref{lem:3.5}, we have
\begin{eqnarray}\label{eq3.2}
\frac{Sz(G)}{W(G)}\geqslant 1+\frac{12(n-3)}{n^3-13n+30}>1+\frac{12}{n^3-19n+54}.
\end{eqnarray}
\begin{figure}[h!]
\begin{center}
\psfrag{a}{$u_3$}\psfrag{b}{$v_3$}\psfrag{c}{$u_4$}
\psfrag{d}{$v_4$}\psfrag{e}{$u_1$}
\psfrag{f}{$v_1$}\psfrag{g}{$u_2$}\psfrag{h}{$v_2$}\psfrag{A}{$(c)$}
\psfrag{B}{$(d)$}\psfrag{C}{$(a)$}\psfrag{D}{$(b)$}\psfrag{x}{$\cdots$}
\includegraphics[width=140mm]{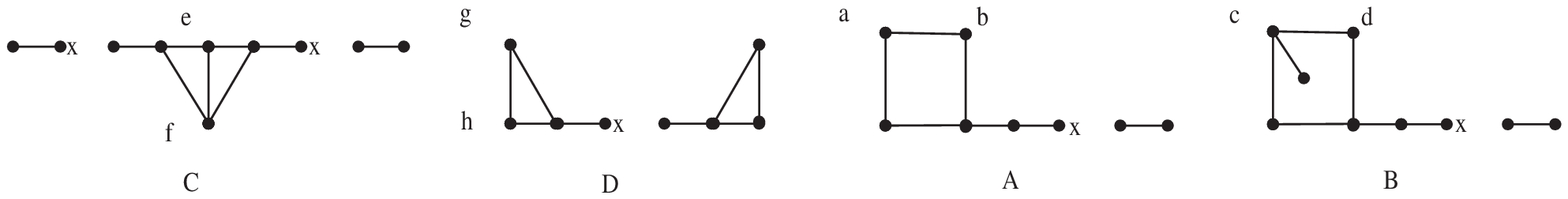}
  \caption{Graphs $(a),(b),(c),(d)$ used in the proofs of Theorems \ref{thm4.1} and \ref{thm4.2}.~}
\end{center}
\end{figure}

If $d(G)\leqslant n-3,$ then $G\ncong B_n^{3,3}$ (based on the fact that $G$ contains at least one non-complete block).
Note that $W(G)<W(G-e)$ for any non-cut edge $e\in E_G.$
Thus, we may repeatedly delete the non-cut edges of $G$ until the resulting graph is a bicyclic graph $G_1$ such that $W(G)\leqslant W(G_1).$ If $G_1\cong B(n,s)$ for some $0\leqslant s\leqslant\lfloor\frac{n-4}{2}\rfloor,$ as $G\ncong  B(n,s),$ then there exist some edges $e_1,e_2,\ldots,e_k,\, k\geqslant 1,$ such that $G=G_1+\{e_1,e_2,\ldots,e_k\},$ where $G_1$ is depicted in Fig.~5(a). Now we construct a new bicyclic graph $G_2$ as follows:
$$
G_2=G_1-u_1v_1+e_1.
$$
Clearly, $G_2$ is not in $\{B_n^{3,3}\}\bigcup \{B(n,s): s=0,1,\ldots, \lfloor\frac{n-4}{2}\rfloor\}.$
It is routine to check that $G=G_2+\{u_1v_1,e_2,\ldots,e_k\}.$ Thus, $W(G)<W(G_2).$

If $G_1\cong B_n^{3,3},$ as $G\ncong B_n^{3,3},$ then there exist some edges $e_1,e_2,\ldots,e_t,\, t\geqslant 1,$ such that $G=G_1+\{e_1,e_2,\ldots,e_t\},$ where $G_1$ is depicted in Fig.~5(b). Now we construct a new bicyclic graph $G_3$ as follows:
$$
G_3=G_1-u_2v_2+e_1.
$$
Clearly, $G_3$ is not in $\{B_n^{3,3}\}\bigcup \{B(n,s): s=0,1,\ldots, \lfloor\frac{n-4}{2}\rfloor\}.$
It is routine to check that $G=G_3+\{u_2v_2,e_2,\ldots,e_t\}.$ Thus, $W(G)<W(G_3).$

Thus, we may assume that $G_1$ is not in $\{B_n^{3,3}\}\bigcup \{B(n,s): s=0,1,\ldots, \lfloor\frac{n-4}{2}\rfloor\}.$ By Theorem \ref{thm5.2}, we have $W(G_1)\leqslant\frac{n^3-19n+54}{6},$
the equality holds if and only if $G_1\cong B_n^{(1)}.$ Together with $W(G)\leqslant W(G_1),$
we can obtain that $W(G)\leqslant\frac{n^3-19n+54}{6}$ with equality if and only if $G\cong B_n^{(1)}.$
By Lemma \ref{lem:2.1}, we have $Sz(G)-W(G)\geqslant 2.$
Therefore,
\begin{eqnarray}
\frac{Sz(G)}{W(G)}-1&=& \frac{Sz(G)-W(G)}{W(G)}\notag\\
                    &\geqslant&\frac{6(Sz(G)-W(G))}{n^3-19n+54} \label{eq3.3}\\
                   &\geqslant &\frac{12}{n^3-19n+54},\label{eq3.4}
\end{eqnarray}
where the equality in (\ref{eq3.3}) holds if and only if $G\cong B_n^{(1)};$ while
the equality in (\ref{eq3.4}) holds if and only if $G\cong G',$ where $G'$ is depicted in Fig. \ref{fig1}.
Hence, $G\cong B_n^{(1)}.$

Therefore, in view of (\ref{eq3.1}), (\ref{eq3.2}) and (\ref{eq3.4}), we have
$\frac{Sz(G)}{W(G)}\geqslant 1+\frac{12}{n^3-19n+54}$
with equality if and only if $G\cong B_n^{(1)},$ as desired.
\qed
\vspace{4mm}

Now, we determine the second smallest value on $Sz^*(G)/W(G)$ among all connected cyclic graphs.
\vspace{2mm}

\noindent{\bf Proof of Theorem \ref{thm4.2}.}\
It is straightforward to check that $\frac{Sz^*(H_n^2)}{W(H_n^2)}=1+\frac{24(n-2)}{n^3-19n+54}.$
Note that $\frac{Sz^*(G)}{W(G)}\geqslant 1+\frac{3(n^2+4n-6)}{2(n^3-7n+12)}$ if $G$ is a non-bipartite graph;
whereas $\frac{Sz^*(G)}{W(G)}\geqslant 1+\frac{24(n-2)}{n^3-13n+36}$ if $G$ is a bipartite graph and the equality holds if and only if $G\cong L_{n,4}.$ Hence, combining with Theorem \ref{lem:2.4}, it suffices to characterize the bipartite graphs with the smallest value on $Sz^*(G)/W(G)$ for $G\ncong L_{n,4}.$ We may complete our proof by comparing this smallest value with $1+\frac{3(n^2+4n-6)}{2(n^3-7n+12)}.$

Let $G$ be a bipartite graph with at least one cycle and $G\ncong L_{n,4}.$ In this case, $Sz^*(G)=Sz(G).$
If $G\cong C_4(P_{n-4},P_1,P_2,P_1),$ then by a direct calculation, we have
$$\frac{Sz^*(G)}{W(G)}=1+\frac{36(n-3)}{n^3-19n+66}>1+\frac{24(n-2)}{n^3-19n+54}.$$

Now, consider $G\ncong C_4(P_{n-4},P_1,P_2,P_1).$ Note that $W(G)<W(G-e)$ for any non-cut edge $e\in E_G.$ Thus, we may repeatedly delete the non-cut edges of $G$ until the resulting graph is a unicyclic graph $G_1$ such that $W(G)\leqslant W(G_1).$ Clearly, $G_1$ is bipartite.
If $G_1\cong L_{n,4},$ as $G\not\cong L_{n,4},$ then there exist some edges $e_1,e_2,\ldots,e_k,\, k\geq1,$ such that $G=G_1+\{e_1,e_2,\ldots,e_k\},$ where $G_1$ is depicted in Fig.~5(c). Then we construct a new unicyclic graph $G_2$ as follows:
$$
G_2=G_1-u_3v_3+e_1.
$$
Clearly, $G_2$ is not in $\{L_{n,4},C_4(P_{n-4},P_1,P_2,P_1)\}.$
It is routine to check that $G=G_2+\{u_3v_3,e_2,\ldots,e_k\}.$ Thus, $W(G)<W(G_2).$

If $G_1\cong C_4(P_{n-4},P_1,P_2,P_1),$ as $G\ncong C_4(P_{n-4},P_1,P_2,P_1),$ then there exist some edges $e_1,e_2,\ldots,\linebreak e_t,\, t\geqslant 1,$ such that $G=G_1+\{e_1,e_2,\ldots,e_t\},$ where $G_1$ is depicted in Fig.~5(d). Then we construct a new unicyclic graph $G_3$ as follows:
$$
G_3=G_1-u_4v_4+e_1.
$$
Clearly, $G_3$ is not in $\{L_{n,4},C_4(P_{n-4},P_1,P_2,P_1)\}.$
It is routine to check that $G=G_3+\{u_4v_4,e_2,\ldots,e_t\}.$ Thus, $W(G)<W(G_3).$

Thus, we may assume that $G_1$ is not in $\{L_{n,4},C_4(P_{n-4},P_1,P_2,P_1)\}.$ By Theorem \ref{thm5.1}, we have $W(G_1)\leqslant \frac{n^3-19n+54}{6},$
the equality holds if and only if $G_1\cong H_n^2,H_n^3$ if $n\geqslant 12$,  $G_1\cong H_{11}^2,\, H_{11}^3$ or $C_3(P_6,P_4,P_1)$ if $n=11$ and $G_1\cong H_{10}^2,$ or $H_{10}^3$ if $n=10$. Together with the fact $W(G)\leqslant W(G_1)$ and $G$ is bipartite,
we can obtain that $W(G)\leqslant\frac{n^3-19n+54}{6}$ with equality if and only if $G\cong H_n^2$ if $n\geqslant 10.$
By Theorem \ref{thm1.1}, we have $Sz^*(G)-W(G)\geqslant 4n-8.$
Therefore,
\begin{eqnarray}
\frac{Sz^*(G)}{W(G)}-1&=& \frac{Sz^*(G)-W(G)}{W(G)}\notag\\
                    &\geqslant &\frac{6(Sz(G)-W(G))}{n^3-19n+54} \label{eq4.5}\\
                    &\geqslant &\frac{24(n-2)}{n^3-19n+54},\label{eq4.6}
\end{eqnarray}
where the equality in (\ref{eq4.5}) holds if and only if $G\cong H_n^2$ if $n\geqslant 10;$ by Theorem \ref{thm1.1}, the equality in (\ref{eq4.6}) holds if and only if $G$ is composed of a cycle $C_4$ on 4 vertices, and one tree rooted at a vertex of the $C_4.$ Hence, $\frac{Sz^*(G)}{W(G)}= 1+\frac{24(n-2)}{n^3-19n+54}$ with equality if and only if $G\cong H_n^2.$

Note that $\frac{3(n^2+4n-6)}{2(n^3-7n+12)}>\frac{24(n-2)}{n^3-19n+54}$ for $n\geqslant 10.$ Thus, we have $\frac{Sz^*(G)}{W(G)}\geqslant 1+\frac{24(n-2)}{n^3-19n+54}$
with equality if and only if $G\cong H_n^2.$

This completes the proof.
\qed

\acknowledgements The authors would like to express their sincere gratitude to all of the referees for their insightful comments and suggestions, which led to a number of improvements to this paper.


\begin{thebibliography}{2}
\providecommand{\natexlab}[1]{#1}
\providecommand{\url}[1]{\texttt{#1}}
\expandafter\ifx\csname urlstyle\endcsname\relax
  \providecommand{\doi}[1]{doi: #1}\else
  \providecommand{\doi}{doi: \begingroup \urlstyle{rm}\Url}\fi


\bibitem[Aouchiche et~al.(2005)Aouchiche, Caporossi, Hansen and Laffay]{1}M. Aouchiche, G. Caporossi, P. Hansen, M. Laffay, AutoGraphiX: a survey, Electron. Notes Discrete Math. 22 (2005) 515-520.
\bibitem[Aouchiche and Hansen (2010)Aouchiche and Hansen]{2}M. Aouchiche, P. Hansen, On a conjecture about the Szeged index, European J. Combin. 31 (2010) 1662-1666.
\bibitem[Bondy and Murty (2008)Bondy and Murty]{3}J.A. Bondy, U.S.R. Murty, Graph Theory, in: GTM, vol. 224, Springer, 2008.
\bibitem[Caporossi and Hansen (2000)Caporossi and Hansen]{4}G. Caporossi, P. Hansen, Variable neighborhood search for extremal graphs 1. The AutoGraphiX system, Discrete Math. 212 (2000) 29-44.
\bibitem[Chen et~al.(2012)Chen, Li, Liu and Gutman]{5}L.L. Chen, X.L. Li, M.M. Liu, I. Gutman, On a relation between the Szeged and the Wiener indices of bipartite graphs, Trans. Comb. 1 (4) (2012) 43-49.
\bibitem[Chen et~al.(2014)Chen, Li and Liu]{C-L-L}L.L. Chen, X.L. Li, M.M. Liu, The (revised) Szeged index and the Wiener index of a non-bipartite graph, European J. Combin. 36 (2014) 237-246.
\bibitem[Cvetkovi\'{c} et~al.(1979)Cvetkovi\'{c}, Boob and Sachs]{C-B-S} D.M. Cvetkovi\'{c}, M. Boob, H. Sachs, Spectral of Graphs. Theory and Application, New York, London, 1979.
\bibitem[Deng (2007)Deng]{D} H. Deng, The trees on $n\geqslant 9$ vertices with the first to seventeenth greatest Wiener
indices are chemical trees, MATCH Commun. Math. Comput. Chem. 57 (2007) 393-402.
\bibitem[Dobrynin et al.(2001)Dobrynin, Entringer and Gutman]{D-E-G} A. Dobrynin, R. Entringer, I. Gutman, Wiener index of trees: theory and applications, Acta Appl. Math., 66 (2001) 211-249.
  \bibitem[Dobrynin and Gutman (1994)Dobrynin and Gutman]{7}A. Dobrynin, I. Gutman, On a graph invariant related to the sum of all distances in a graph, Publ. Inst. Math.(Beograd) (N.S.) 56 (1994) 18-22.
\bibitem[Dobrynin and Gutman (1995)Dobrynin and Gutman]{8}A. Dobrynin, I. Gutman, Solving a problem connected with distances in graphs, Graph Theory Notes N. Y. 28 (1995) 21-23.
\bibitem[Dobrynin et al.(2002)Dobrynin,Gutman,Klav\v zar,\v Zigert]{D-G-K-Z} A.A.~Dobrymin, I.~Gutman, S.~Klav\v zar, P.~\v Zigert, Wiener index of hexagonal systems, Acta Appl. Math. 72 (2002) 247-294.
\bibitem[Dong and Zhou (2012)Dong, Zhou]{D-Z}H. Dong, B. Zhou, Maximum Wiener index of unicyclic graphs with fixed maximum degree, Ars Combin. 103 (2012) 407-416.
 \bibitem[Du and Ili\'c(2013)Du, Ili\'c(2013)]{9}Z. Du, A. Ili\'c, On AGC conjectures regarding average eccentricity, MATCH Commun. Math. Comput. Chem. 69 (2013) 579-587.
\bibitem[Entringer et al.(1976)Entringer, Jackson, Snyder]{008}R.C.~Entringer, D.E.~Jackson, D.A.~Snyder, Distance in graphs, Czech.\ Math.\ J.\ 26 (1976) 283--296.
\bibitem[Gutman (1994)Gutman]{10}I. Gutman, A formula for the Wiener number of trees and its extension to graphs containing cycles, Graph Theory Notes N. Y. 27 (1994) 9-15.
\bibitem[Gutman and Polansky (1986)Gutman, Polansky]{12}I. Gutman, O.E. Polansky, Mathematical Concepts in Organic Chemistry, Springer, Berlin, 1986.
\bibitem[Hansen et al.(2010)Hansen]{Han}P. Hansen et al. Computers and conjectures in chemical graph theory, in: Plenary Talk in the International Conference on Mathematical Chemistry, August 4-7, 2010, Xiamen, China.
\bibitem[Ili\'c (2010)Ili\'c]{15}A. Ili\'c, Note on PI and Szeged indices, Math. Comput. Modelling 52 (2010) 1570-1576.
\bibitem[Khodashenas et al.(2011)Khodashenas, Nadjafi-Arani, Ashrafi and Gutman]{16}H. Khodashenas, M.J. Nadjafi-Arani, A.R. Ashrafi, I. Gutman, A new proof of the Szeged-Wiener theorem, Kragujevac J. Math. 35 (2011) 165-172.
\bibitem[Klav\v{z}ar and Nadjafi-Arani (2013)Klav\v{z}ar, Nadjafi-Arani]{17}S. Klav\v{z}ar, M.J. Nadjafi-Arani, Wiener index versus Szeged index in networks, Discrete Appl. Math. 161 (2013) 1150-1153.
\bibitem[Klav\v{z}ar and Nadjafi-Arani (2014)Klav\v{z}ar, Nadjafi-Arani]{28}S. Klav\v{z}ar, M.J. Nadjafi-Arani, Improved bounds on the difference between the Szeged index and the Wiener index of graphs, European J. Combin. 39 (2014) 148-156.
\bibitem[Klav\v zar et al.(1996)Klav\v zar, Rajapakse and Gutman]{klavzar-1996}S.~Klav\v zar, A.~Rajapakse, I.~Gutman, The Szeged and the Wiener index of graphs, Appl.\ Math.\ Lett.\ 9 (1996) 45--49.
\bibitem[Knor et al.(2016)Knor, \v{S}krekovski and Tepeh]{KM-S}M. Knor, R. \v{S}krekovski, A. Tepeh, Mathematical aspects of Wiener index, Ars Math. Contemp. 11 (2016) 327-352.
\bibitem[Li and Song (2014)Li, Song]{LS-S}S.C. Li, Y.B. Song, On the sum of all distances in bipartite graphs, Discrete Appl. Math.  169  (2014) 176-185.
\bibitem[Li and Zhang (2017) Li, Zhang]{L-Z} S.C. Li, H.H. Zhang, Proofs of three conjectures on the quotients of the (revised) Szeged index and the Wiener index and beyond, Discrete Math. 340 (2017) 311-324.
\bibitem[Li and Liu (2013)Li, Liu]{18} X.L. Li, M.M. Liu, Bicyclic graphs with maximal revised Szeged index, Discrete Appl. Math. 161 (2013) 2527-2531.
\bibitem[Liu et al. (2016)Liu, Du, Jia]{LR-D}R.F. Liu, X. Du, H.C. Jia, Wiener index on traceable and Hamiltonian graphs, Bull. Aust. Math. Soc. 94 (2016) 362-372.
\bibitem[McKay and  Piperno(2014) McKay, Piperno]{M-V0}B.D. McKay, A. Piperno, Practical Graph Isomorphism, II,
  J. Symbolic Comput. 60 (2014) 94-112, \url{http://dx.doi.org/10.1016/j.jsc.2013.09.003}
\bibitem[Mukwembi and Vetr\'{i}k (2014)Mukwembi, Vetr\'{i}k]{M-V} S. Mukwembi, T. Vetr\'{i}k, Wiener index of trees of given order and diameter at most 6, Bull. Austra. Math. Soc. 89 (2014) 379-396.
\bibitem[Nadjafi-Arani et al. (2011)Nadjafi-Arani, Khodashenas and Ashrafi]{19}M.J. Nadjafi-Arani, H. Khodashenas, A.R. Ashrafi, On the differences between Szeged and Wiener indices of graphs, Discrete Math. 311 (2011) 2233-2237.
\bibitem[Nadjafi-Arani et al. (2012)Nadjafi-Arani, Khodashenas and Ashrafi]{20}M.J. Nadjafi-Arani, H. Khodashenas, A.R. Ashrafi, Graphs whose Szeged and Wiener numbers differ by 4 and 5, Math. Comput. Modelling 55 (2012) 1644-1648.
\bibitem[Pisanski and Randi\'{c} (2010)Pisanski and Randi\'{c}]{21}T. Pisanski, M. Randi\'{c}, Use of the Szeged index and the revised Szeged index for measuring network bipartivity, Discrete Appl. Math. 158 (2010) 1936-1944.
\bibitem[Randi\'{c} (2002)Randi\'{c}]{23}M. Randi\'{c}, On generalization of Wiener index for cyclic structures, Acta Chim. Slov. 49 (2002) 483-496.

\bibitem[Simi\'{c} et al.(2000)Simi\'{c}, Gutman and Balti\'{c}]{24}S. Simi\'{c}, I. Gutman, V. Balti\'{c}, Some graphs with extremal Szeged index, Math. Slovaca 50 (2000) 1-15.
\bibitem[Wiener (1947)Wiener]{25}H. Wiener, Structural determination of paraffin boiling points, J. Am. Chem. Soc. 69 (1947) 17-20.
\bibitem[Xing and Zhou (2011)Xing and Zhou]{26}R. Xing, B. Zhou, On the revised Szeged index, Discrete Appl. Math. 159 (2011) 69-78.
\bibitem[Xu and Das (2013)Xu and Das]{X-D}K.X. Xu, K.C. Das, Extremal unicyclic and bicyclic graphs with respect to Harary index, Bull. Malays. Math. Sci. Soc. 36 (2013), 373-383.
\bibitem[Tang and Deng(2008) Tang and Deng]{T-D} Z.K. Tang, H.Y. Deng, The $(n,n)$-graphs with the first three extremal Wiener indices, J. Math. Chem. 43 (2008) 60-74.
\bibitem[Zhang et al. (2016)Zhang, Li and Zhao]{Z-L-Z} H.H. Zhang, S.C. Li, L.F. Zhao, On the further relation between the (revised) Szeged index and the Wiener index of graphs, Discrete Appl. Math. 206  (2016) 152-164.
\bibitem[Zhang et al.(2010)Zhang, Liu and Han]{Z} X.D. Zhang, Y. Liu, M.X. Han, Maximum Wiener index of  trees with given degree sequence,
  MATCH Commun. Math. Comput. Chem. 64 (2010) 661-682.

\end{thebibliography}
\end{document}